\documentclass [12pt,reqno]{article}
\usepackage[leqno]{amsmath}
\usepackage{paralist, amsthm,amsfonts,amsmath,amssymb,graphics,hyperref, pstricks,graphicx
,
multicol,
}%
\usepackage{t1enc, pst-poly,pstricks}

\divide
\oddsidemargin by 2 
\makeatletter
\newcommand{\address}[2][]{%
  \ifx\@add@ress\@undefined\gdef\@add@ress{\par\par\bigskip}\AtEndDocument{\@add@ress}\fi
  \g@addto@macro\@add@ress{\bigskip\noindent{\small\scshape%
      \ifx#1\empty\else{\bfseries Address of #1:}\ \fi#2}\par\par}}

\renewenvironment{abstract}{\small\quotation\noindent
  {\bfseries \abstractname}}{\endquotation \par}
\newcommand{\footnotetextplain}[1]{\begingroup\def\@thefnmark{}%
  \long\def\@makefntext##1{\parindent 0pt\noindent ##1}\@footnotetext{#1}
  \endgroup}

\newcommand{\AMSsubjclass}[2]{\footnotetextplain{2000
   \emph{Mathematics Subject Classification:} Primary #1, Secondary #2.}}

\newcommand{\keywords}[1]{\footnotetextplain{\emph{Key words and phrases:} #1.}}

\makeatletter \xdef\qedbuit{\qed}

\newcommand{\TeoremaAmbFinalMarcat}[1]{%
  \expandafter\gdef\csname end#1\endcsname{\qedbuit\@endtheorem}}
\newcommand{\vs}{\vspace{.1in}}  
\newenvironment{prooftext}[1] 
               {\noindent \textit{ #1}~ }     
               {\hfill\rule{2.5mm}{2.5mm} \vspace{\parskip} } 

\theoremstyle{definition}
 \TeoremaAmbFinalMarcat{defi}
 \TeoremaAmbFinalMarcat{ex}
 \TeoremaAmbFinalMarcat{rem}
 \TeoremaAmbFinalMarcat{conv}

\newenvironment{proclama}[1]{\par\vspace{\topsep}\noindent{\bf #1}
 \begin{em}}{\end{em}\par\vspace{\topsep}}

\newcommand{\start}[2]{\begin{#1}\label{#2}}

\newcommand{\secc}[1]{Section~\ref{#1}}
\newcommand{\theoc}[1]{Theorem~\ref{#1}}
\newcommand{\propc}[1]{Proposition~\ref{#1}}
\newcommand{\coryc}[1]{Corollary~\ref{#1}}
\newcommand{\lemc}[1]{Lemma~\ref{#1}}
\newcommand{\defic}[1]{Definition~\ref{#1}}
\newcommand{\remc}[1]{Remark~\ref{#1}}

\newcommand{\figc}[1]{Figure~\ref{#1}}
\newcommand{\tabc}[1]{Table~\ref{#1}}

\newcommand{\ctc}[1]{Computational Theorem~\ref{#1}}

\def\@enum@{\list{\csname label\@enumctr\endcsname}%
           {\usecounter{\@enumctr}\def\makelabel##1{\hss\llap{##1}}
           \itemsep=2pt\parsep=0pt\topsep=3pt plus 1pt minus 1 pt}}
\newenvironment{alphlist}{\enumerate[(a)]}{\endenumerate}
\newenvironment{romlist}{\enumerate[(i)]}{\endenumerate}
\newenvironment{numlist}{\enumerate[(1)]}{\endenumerate}
\newenvironment{lplist}{ \enumerate[I.]}{\endenumerate}

\newcommand{\si}{\mathop\mathrm{SI}}
\newcommand{\inn}{\mathrm{IN}}
\newcommand{\A}{A}
\newcommand{\B}{B}

\newcommand{\uu}{\mathsf u}
\newcommand{\vv}{\mathsf v}
\newcommand{\BB}{\mathsf B}
\newcommand{\AfA}{\mathsf A}
\newcommand{\bb}{\mathsf b}
\newcommand{\afa}{\mathsf a}

\newcommand{\VV}{\mathsf V}
\newcommand{\YY}{\mathsf Y}
\newcommand{\rr}{\mathsf r}
\newcommand{\RR}{\mathsf R}

\newcommand{\xx}{\mathsf x}
\newcommand{\XX}{\mathsf X}
\newcommand{\zz}{\mathsf z}
\newcommand{\ZZ}{\mathsf Z}
\newcommand{\ww}{\mathsf w}
\newcommand{\WW}{\mathsf X}

\newcommand\yy{\mathsf y}

\newcommand{\ra}{\rangle}
\newcommand{\la}{\langle}
\newcommand{\bnl}{\begin{numlist}}
\newcommand{\enl}{\end{numlist}}

\title{Self-intersection numbers of curves\\ on the
punctured torus}
\author{Moira Chas and Anthony Phillips}
\date{\today}

\begin{document}
\maketitle

\begin{abstract} 
The minimum number of self-intersection points for members of a free homotopy class of curves on the punctured torus  is bounded above in terms of the number $L$ of letters required for a minimal description of the class in terms of the generators of the fundamental group and their inverses: it is less than or equal to $(L-2)^{2}/4$ if L is even, and $(L-1)(L-3)/4$ if L is odd. The classes attaining this bound are explicitly described in terms of the generators; there are $(L-2)^2+ 4$ of them if L is  even, and $2(L-1)(L-3)+8$ if $L$ is odd; similar descriptions and totals are given   for classes with self-intersection number equal to one less than the maximum. 

Proofs use both combinatorial calculations and  topological operations on representative curves. Computer-generated data are tabulated counting, for each
non-negative integer, how many length-$L$ classes 
have that self-intersection number, for each length $L$ less than or 
equal to 12. 
Experimental data are also presented for the pair-of-pants surface. 

\end{abstract}

\AMSsubjclass{57M05}{57N50,30F99}

\section{Introduction}

The punctured torus has the homotopy type of a figure-eight.
Its fundamental group is free on two generators: once 
these are chosen, say  $a, b$, 
a free homotopy class of curves on the surface can
be uniquely 
represented as a reduced cyclic word in the symbols
$a, b, \A, \B$ (where $\A$ stands for $a^{-1}$ and $\B$ for $b^{-1}$). 
A {\em cyclic word } $w$ is an equivalence class of words 
related by a cyclic permutation of their letters; we will write
$w=\langle r_1 r_2 \dots r_n \rangle$ where the $r_i$ are
the letters of the word, and $\langle r_1 r_2 \dots r_n \rangle =
\langle r_2 \dots r_n r_1 \rangle$, etc. {\em Reduced}
means that the cyclic word contains no 
juxtapositions of $a$ with $A$, or $b$ with $B$.
Note here that
we will call a free homotopy class (a reduced cyclic word) {\em primitive} 
if is not a proper power of
another class (another word); and among the non-primitive classes are
words we will call {\em pure powers:} those which are a proper power
of a generator.
The {\em length} (with respect to the generating set $(a,b)$) of 
a free homotopy class of curves is the number of letters occurring in
the corresponding reduced cyclic word. 


This work studies the relation between length and
the {\em self-intersection
number} of a free homotopy class of curves: the 
smallest number of 
self-intersections among all general-position curves in the class. 
(General position in this context means as usual that there are
no tangencies or multiple intersections).
The self-intersection number is a property of the free homotopy class
and hence 
of the corresponding reduced cyclic word $w$; 
we denote it by $\si (w)$.

\start{theo}{upper bound 2}

\begin{figure}[http]\label{grids}
\centerline{\includegraphics[width=5in]{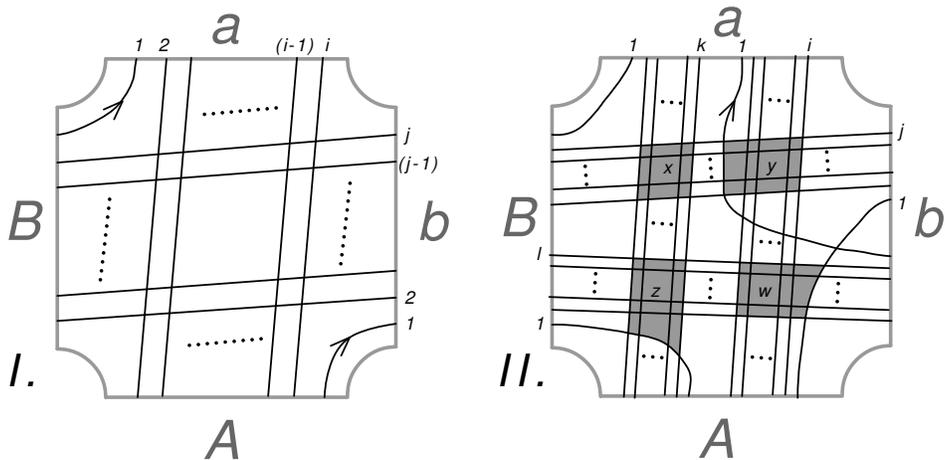}}

\caption{Curves of maximal self-intersection on the punctured torus.
I. $w = \la a^ib^j\ra $ with $(i-1)(j-1)$ intersection points, a maximum
when $i=j$ (even length) or $i=j\pm 1$ (odd length). II.
$w = \la a^ib^ja^kB^l\ra $. Block ``x'' has $(k-1)(j-1)$ intersection points;
block ``y'' has $i(j-1)$; block ``z'' has $l(k-1)$; block ``w''
has $i(l-1)$; and there are an additional $i$. The total is
$(i+k-1)(j+l-1)$, a maximum when $i+k=j+l$ (even length) or
$i+k = j+l \pm 1$ (odd length). Graphic conventions from
\secc{prelim}; curve II drawn using the algorithm
of \cite{blood}. Similar diagrams appear in \cite{chemotti}. }
\end{figure}

The maximal self-intersection
number for a primitive reduced cyclic word of 
length $L$ on the punctured torus is:
$$\left \{ \begin{array}{ll} (L-2)^2/4 & \mbox{if $L$ is even,} \\
(L-1)(L-3)/ 4 &   \mbox{if $L$ is odd.}
\end{array} \right . $$ 
The words realizing the maximal self-intersection number are 
(see \figc{grids}): \begin{numlist}
\item[{\rm (1)}] $L$ even: \begin{romlist}
    \item[{\rm (i)}] $\la r^{L/2}s^{L/2} \ra$, $r \in \{a, A\},~ s\in \{b, B\}$
    \item[{\rm (ii)}] $\la r^is^jr^{L/2-i}S^{L/2-j} \ra$, $r \in \{a, A\},~ s\in \{b, B\}$,~$S=s^{-1}$, and similar configurations
interchanging $r$ and $s$.
\end{romlist}
\item[{\rm (2)}] $L$ odd: \begin{romlist}
\item[{\rm (i)}] $ \la r^{(L+1)/2}s^{(L-1)/2} \ra$, $r \in \{a, A\},~ s\in \{b, B\}$, or vice-versa
\item[{\rm (ii)}] $\la r^is^jr^{(L+1)/2-i}S^{(L-1)/2-j} \ra$,
$\la r^is^jr^{(L-1)/2-i}S^{(L+1)/2-j} \ra$,
$r \in \{a, A\}, ~s\in \{b, B\}$,~ $S=s^{-1}$, and similar configurations
interchanging $r$ and $s$.
\end{romlist}
\end{numlist}
\end{theo}

Elementary counting with Theorem \ref{upper bound 2} yields the next result:

\start{theo}{Counting Theorem}
The number of distinct primitive 
free homotopy classes of length $L$ realizing the
maximal self-intersection number is 
$$\left \{ \begin{array}{ll} (L-2)^2+4 & \mbox{if $L$ is even,} \\
2(L-1)(L-3)+8&   \mbox{if $L$ is odd.}
\end{array} \right . $$ \label{counting}
\end{theo}


Elementary computation with Theorem \ref{counting} allows the
inequality to be reversed:
 
\start{theo}{lblt}
Let $w$ be the 
 reduced cyclic word corresponding to a primitive
free homotopy class of curves on the punctured torus. Then if
$\si(w)\geq 1$,
the length of $w$ is greater than or equal to the smallest integer 
larger than  $2\sqrt{\si(w)}+2$. Moreover, this bound is sharp.
\end{theo}

\start{rem}{rem1}
Pure-power words of length $L$ between 2 and 6 do
not fit the pattern of Theorems \ref{upper bound 2} -  \ref{lblt}. 
Namely, $\si(r^L)=L-1 > (L-2)^2/4 \mbox{~and~} (L-1)(L-3)/4$ for
integers in that range. Theorems \ref{upper bound 2} - \ref{lblt}  can in fact be extended to all words of length seven or more,  primitive
or not.
\end{rem}

\start{rem}{generating} The length of the word representing a free homotopy class depends
on the choice of generating set $(a, b)$ for the fundamental group,
while its self-intersection number does not. Since the theorems
above apply for {\em any} generating set, they can be rephrased
in terms of  the shortest of such lengths. 
\end{rem}

\start{rem}{generating1} The group of automorophisms of the fundamental group of the punctured torus acts on the set of cyclic words with a fixed self-intersection number $n$. Words with maximal self-intersection number minimize length in an orbit of this action. Rivin asked us if every orbit contains a word with maximal self-intersection number for its length.  But
$w=\la ababAB \ra$ is not in the orbit of such a word (this can be proved using \cite[Proposition 4.19]{LS}).
\end{rem}

Theorems \ref{upper bound 3} and \ref{counting 2} treat curves,
on the punctured torus,
of self-intersection number one less than the maximum for their
length; we do not have similar formulas for
the distribution of other self-intersection numbers 
among curves of a given length. Here is some numerical
evidence, computed using 
the algorithm given in \cite{cl}; see also \cite{chas}; the Java program can be found at \cite{chas program}. This evidence was in fact the motivation for the research
presented here.

\start{ct}{comp}
The number of distinct
primitive free homotopy classes 
with a given number of
self-intersections, corresponding to primitive reduced cyclic words of a given length appears, for length up to  $12$, in 
\tabc{counting table}. (If one entry of a row of \tabc{counting table} is 
$0$ then all the entries to its right are also $0$.)
\end{ct}

\begin{table}[htbp]
   \centering
  \begin{tabular}{|c|r|r|r|r |r|r|r|r|r |r|r|}
\hline
$\mathrm{length}\setminus  \mathrm{SI}$&0&1&2&3&4&5&6&7&8&9&10\\\hline
1&{\bf 4}&0&0&0&0&0&0&0&0&0 &0\\\hline 
2&{\bf 4}&0&0&0 &0&0&0&0&0&0&0 \\\hline 
3&{\bf 8}&0&0&0 &0&0&0&0&0&0&0 \\\hline 
4&{\it 10}&{\bf 8}&0&0 &0&0&0&0&0&0&0 \\\hline 
5&16&{\it 8}&{\bf 24}&0 &0&0&0&0&0 &0&0\\\hline 
6&8&16&32&{\it 40 }&{\bf 20}&0&0&0&0 &0&0\\\hline 
7&24&16&32&48 &112&{\it 24}&{\bf 56}&0&0 &0&0\\\hline 
8&16&24&52&76 &116&156&136&104&{\it 90 }&{\bf 40}&0\\\hline 
9&24&32&64&120 &144&240&384&208&376 &136&304\\\hline 
10&16&32&72&168 &272&332&492&628&644 &700&700\\\hline 
11&40&48&80&160 &272&584&664&1200&1280 &1368&1608\\\hline 
12&16&40&104&208&372&660&1048&1408&2044&2696&3088\\\hline 
\end{tabular}\vspace{0.5cm}
\begin{tabular}{|c|r|r|r|r|r|r|r|r|r|}\hline
$\mathrm{length} \setminus \mathrm{SI}$ &11&12&13&14&15& 16&17&18&19 \\\hline
9&{\it 48}&{\bf 104}&0 &0&0&0&0&0 &0\\\hline
10&548&464&360 &224&{\it 160}&{\bf 68}&0&0&0 \\\hline
11&1368&2048&976 &1704&528&1072&264&  592 &{\it 80}\\\hline
12&3580&3866&3792&3816&3612&3272&2820&2276&1808\\\hline
\end{tabular}\vspace{0.5cm}
\begin{tabular}{|c|r|r|r|r|r|r|r|}\hline
$\mathrm{length} \setminus \mathrm{SI}$&20&21 &22&23&24&25&26\\\hline
11&{\bf 168}&0&0&0&0&0&0\\\hline
12 &1308&960& 680&392&{\it 250}&{\bf 104}&0\\\hline
\end{tabular}

\caption{The $i,j$ entry in this table is the number of
distinct reduced primitive cyclic words of length $i$ with
exactly $j$ self-intersections, up to the maximum possible
self-intersection number for each length. {\bf Bold-face} numbers
and their location 
correspond to Theorems  \ref{upper bound 2} 
and \ref{counting}, {\em italic} numbers to \theoc{counting 2}.}
\label{counting table}
\end{table}

\start{ct}{cct}
Let  $k \in \{1,2,\dots,25\}$ and let $K$ be  
the set of all cyclic reduced words $v$ corresponding to primitive
free homotopy classes of curves on the punctured torus, with  $\si(v) \ge k$. 
If $w$  is a word in $K$ 
with minimal length then the following statements hold:
\begin{numlist} 
\item[{\rm(1)}] The length of $w$ is equal to the smallest integer larger than or equal to $2\sqrt{k}+2$. 
\item[{\rm(2)}] $\si(w)=k$.
\end{numlist}
\end{ct}


\subsection{Related results}\label{rr}

For a 
reduced cyclic word $w$ written in the symbols $\{a, A, b, B\}$, let $\alpha(w)$ and $\beta(w)$ denote
the total number of occurrences of $a,A$ and of $b,B$, respectively. 
Andrew Blood \cite{blood} gives a simple
construction of a representative curve which has  at most 
$(\alpha(w)-1)(\beta(w)-1)$ intersections; he also finds some of the words 
whose representative curves
require this number of self-intersections, 
namely those of the form $a^{\alpha}b^{\beta}$. Together these
two discoveries constitute a different proof of the first part of 
our \theoc{upper bound 2} (compare \theoc{main-th}). In addition,
Frank Chemotti and Andrea Rau \cite{chemotti} give elementary proofs of parts
(2), (3) and (4) of our \propc{lp}.
 This unpublished
work only came to our attention during the final editing of this paper.

Birman and Series \cite{BS}  give an algorithm to decide whether a simple representative exists for a
reduced cyclic word in the generators of the fundamental group of a surface with boundary. These ideas are extended by Cohen and Lustig \cite{cl} (see also \cite{chas} and \cite{Tan}), who give an algorithm to compute the self-intersection of a reduced cyclic  word. The program to compute \tabc{comp} is based on these algorithms.

From the geometric point of view, the punctured torus has been studied as 
a manifold with boundary:
the complement in $S^1\times S^1$ of an open disc. 
This manifold admits a complete hyperbolic metric for
which the boundary circle is a geodesic. Since every
free homotopy class contains exactly one geodesic representative, 
and since a primitive
geodesic cannot have excess intersections \cite{hs}, 
the results in this section translate into
results about counting geodesics on that Riemann surface.

\begin{enumerate}

\item[(1)] It follows from Cohen and Lustig \cite[Main Theorem]{cl} (see also \cite[Proposition 2.9 and Remark 3.10]{chas})
that for any surface $S$ with non-empty boundary and  negative
Euler characteristic, $\si(w) \leq L(L-1)/2$ (using our notation)
for $w$ a primitive word of length $L$ in the generators (and their inverses) of the fundamental group of $S$. For the torus with one boundary component,
the special case examined here, our upper bound
(Theorem \ref{counting}) is lower.
\end{enumerate}

For the torus with one geodesic boundary, 
once a  pair of free generators is chosen
for  the fundamental group 
then any hyperbolic metric, restricted to closed geodesics, 
is quasi-isometric to the word-length metric. This is a special case of the
\v{S}varc-Milnor Lemma \cite{Milnor}, \cite{Bridson}.
Hence we can refer to word-length as  {\em combinatorial length}.

\begin{enumerate}
\item[(2)] Lalley \cite[Theorem 1]{lalley} proved that on a compact, 
hyperbolic, closed 
surface {\em most}  closed geodesics of length  
approximately $\ell$ have about $C \ell^{2}$ self-intersections for some 
positive constant $C$ depending on the surface. 
As a consequence of 
our \theoc{upper bound 2}, in the case of the torus with one geodesic boundary
component, for each hyperbolic metric there exists a positive constant $C'$ 
such that the number of self-intersection points of {\em every} geodesic 
of length $\ell$ is less than $C' \ell^{2}$. (This fact also admits
an elementary proof, as Lalley pointed out to us). Lalley 
also studies the distribution on the surface of self-intersection
points of a typical geodesic; \cite[Theorem 2]{lalley} 
may be compared with the patterns in \figc{grids}.

\item[(3)] Basmajian proves in \cite[Corollary 1.2]{basmajian} that for any hyperbolic surface 
 there exists an increasing sequence of 
constants $\{M_k\},  k \ge 1$, 
tending to infinity so that if $\omega$ is a closed   
geodesic with self-intersection number $k$, then the hyperbolic
length  of $\omega$ is greater than $M_{k}$.
For the punctured torus and combinatorial length our \coryc{lblt} 
gives explicit values
for $M_k$, and our bounds are sharp.
\end{enumerate}

In view of the quasi-isometry between combinatorial and hyperbolic length 
for the torus with one boundary component 
 the numbers in \ctc{comp}
are concordant with numbers
or estimates from several other lines of research:
\begin{numlist}

\item[{\rm (1)}] It is known that for any hyperbolic surface the total number of 
primitive closed geodesics of length at most $L$ 
is asymptotic to $e^{hL}/L$ ($h$ is the topological entropy
of the geodesic flow; see \cite{buser} and references therein;
similar results hold for the variable curvature case, \cite{lalley2}, \cite{marg}, \cite{pp}).  On the punctured torus,  the number of 
distinct primitive classes of combinatorial length $L$ at most twelve, i.e.
the sum of the numbers in row $L$ of \tabc{counting table}, appears to be 
very rapidly asymptotic 
to $3^L/L$. 

\item[{\rm (2)}]  The numbers in the first column of \tabc{counting table}, giving the number of simple classes for a given length, can be
compared with the results of McShane and Rivin \cite{rivin mcshane} for the  punctured  torus and Mirzakhani \cite{mm} for a general surface of negative Euler characteristic (see also \cite{rivin} for historical brackground). Mirzakhani, McShane  and Rivin prove that the number of simple closed geodesics of hyperbolic length at most $L$ grows as a quadratic polynomial in $L$ (contrast with
\theoc{Counting Theorem}, where the number of maximal curves of length {\em exactly} $L$ grows quadratically with $L$).
For the range of \tabc{counting table}, we have data consistent with these:  the number of simple curves of length exactly
$2n+1$, $n\geq 1$, appears to grow more or less linearly with $n$; for $2n+1$ a prime, it is exactly $8n$. 

\item[{\rm (3)}]  For $L$ even, the numbers in the second column of \tabc{counting table} grow as $4(L-2)$. This is consistent with Rivin's \cite{rivin single} determination that the number of single-self-intersection geodesics of length
at most $L$ grows quadratically with $L$.

\item[{\rm (4)}] For a closed surface $S$, Basmajian proves in \cite[Proposition 1.3]{basmajian}
 that there are constants $N_{k}$ (depending on the genus of 
$S$) such that the shortest geodesic on $S$ with at least $k$ 
intersection points has length bounded above by $N_{k}$.
This generalizes Buser's proof  \cite{buser} that $N_1 = 1$.  
\ctc{cct} gives $N_k$ an explicit value  
for curves of
combinatorial length less than $13$  on 
the punctured torus.

\end{numlist}


\subsection{Sketch of proof}\label{sketch}

The method of proof in this paper keeps track of three integer
parameters of a reduced cyclic word $w$ in the alphabet
$a, b, A, B$: along with
$\alpha(w)$ and  $\beta(w)$ (see \secc{rr}) 
there is
$h(w)$, the total number of block-pairs in $w$;
these are defined as follows:

\start{defi}{Definition 3: block-pairs.}
A reduced cyclic word  $w$ is either a pure power
or there exist pairs of positive
integers $j_1, k_1, \dots j_n, k_n$, $n\ge 1$, such that
  $ w = \langle r_{1}^{j_{1}}s_{1}^{k_{1}} r_{2}^{j_{2}}s_{2}^{k_{2}}
\dots  r_{n}^{j_{n}}s_{n}^{k_{n}}\rangle$,
where $r \in \{a,\A\}$ and $s \in \{b,\B
\}$. Each of the $r_{i}^{j_{i}}s_{i}^{k_{i}}$ occurring in this
expression is a {\em block-pair}; the
{\em number of block-pairs} of $w$  is defined to be $n$ in the second case,
and zero in the first.\label{def-block-pair}
\end{defi}

The main theorem in this paper is \theoc{main-th}; it will be proved in \secc{proof-ub}.

\start{theo}{main-th}  For the punctured torus, let $w$ be the reduced cyclic word 
corresponding to a
free homotopy class of curves with a positive number $h$ of block-pairs.
If $h=1$ then $\si(w)= (\alpha(w)-1)(\beta(w)-1)$. If $h \geq 2$,
 $$\si(w)\le (\alpha(w)-1)(\beta(w)-1)-h+2.$$

The words $w$ realizing the maximal self-intersection for
non-pure-power words with given $\alpha$ and $\beta$ (that is, 
$\si(w)= (\alpha(w)-1)(\beta(w)-1)$) have one of the following forms,
\begin{numlist}
\item[{\rm (1)}] $\la r^is^j \ra$, $r \in \{a, A\}, s\in \{b, B\}$; here $\alpha(w) = i>0,
\beta(w) = j>0$. 
\item[{\rm (2)}] $\la r^i s^j r^k S^l \ra$, all $i,j,k,l>0$,  
where $r \in \{a, A\}$ (and then $i+k=\alpha(w)$)
and $s \in \{b, B\}$ (and then $j+l=\beta(w)$), 
or vice-versa.
\end{numlist}
\end{theo}

This theorem has two immediate corollaries:

\start{cory}{upper bound}
Let $w$ be the reduced cyclic word corresponding to a primitive
free homotopy class of curves on the punctured torus.  Then
$$\si(w)\le (\alpha(w)-1)(\beta(w)-1).$$
\end{cory}

\start{cory}{max-form} Among primitive words those of maximal self-intersection number for their $\alpha$
and $\beta$ values, i.e. with $\si(w)= (\alpha(w)-1)(\beta(w)-1)$,
 have one of the following forms:
\begin{numlist}
\item[{\rm (1)}] $\la r \ra$, $r \in \{a,b,A,B\}$,
\item[{\rm (2)}] $\la r^is^j\ra,$ $r \in \{a, A\}, s\in \{b, B\}$; here $\alpha(w) = i>0,
~\beta(w) = j>0$. 
\item[{\rm (3)}] $\la r^i s^j r^k S^l \ra$, all $i,j,k,l>0$,  
where $r \in \{a, A\}$ (and then $i+k=\alpha(w)$)
and $s \in \{b, B\}$ (and then $j+l=\beta(w)$), 
or vice-versa.
\end{numlist}
\end{cory}

\start{rem}{x1} Since $\alpha(w) + \beta(w) = L$, the length
of $w$, an elementary calculation leads from 
\coryc{upper bound} and \coryc{max-form} to \theoc{upper bound 2}.
\end{rem}

The next three sections carry the proof of \theoc{main-th}.
The strategy is to show that only words of the types listed in the
statement of the theorem, i.e.
 $\la r^is^j \ra$ and $\la r^i s^j r^k S^l \ra$,~ $r \in \{a, A\},
~ s\in \{b, B\}$, or vice-versa,  can have maximum self-intersection number for their length; this 
will be done by exhibiting, for any word which is not of these types,
another word of the same length and with strictly larger 
self-intersection number. For most
words $w$, 
``cross-corner surgery'' (defined below) will produce a $w'$
with the same $\alpha$ and $\beta$ values (and so of the same length),
with $\si(w')>\si(w)$ and with $h(w')<h(w)$
For certain words with two,
three or four blocks, not candidates for surgery, the self-intersection
number will be computed explicitly by counting ``linked pairs'' of
subwords (definition below) and determining that it is indeed smaller
than the self-intersection number of a word of the same length but
of one of the two listed types (whose
self-intersection numbers are also computed
by counting linked pairs).

This work benefited from discussions with Ara Basmajian, 
Joel Hass, Stephen Lalley, Igor Rivin and Dennis Sullivan.

\section{Cross-corner surgery}\label{sec-[httb]}

\subsection{Preliminaries}
\label{prelim}
\begin{figure}[htp]
\centerline{\includegraphics[width=4in]{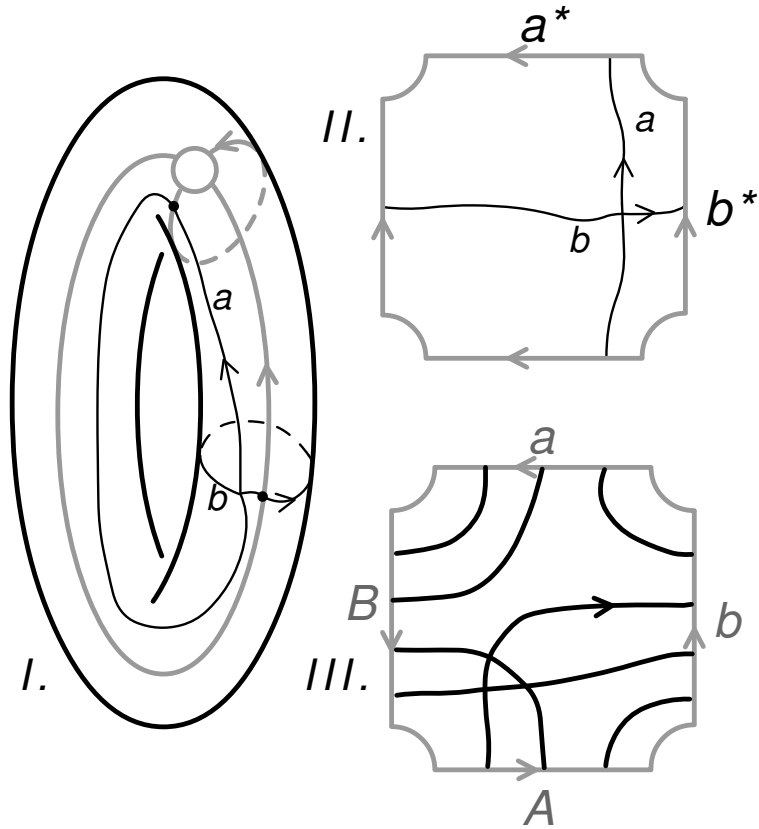}}
\caption{ The punctured torus as a polygon with identifications.
I, II. The generators $a,b$ and their inverses $A,B$ can be identified by their
intersections with the dual cycles $a*, b*$, which appear among the edges
of the fundamental polygon.
III. When an oriented curve
has been lifted to the fundamental polygon,
the cyclic word corresponding to its free homotopy class can
can be obained by choosing a starting point and recording
in sequence the edges it crosses, reading their names from
{\em inside} the polygon. The lifted curve $\langle baBBAba\rangle$
with self-intersection number 3, is shown as an example.}
\label{identifications}
\end{figure}

Here, let $M$ represent the punctured torus as a topological space.
The choice of generators $(a,b)$ for $\pi_1 M$ naturally implies
a {\em fundamental polygon} from which $M$ may be reconstructed by
edge-identification. Namely,
we can choose, as representative cycles for
 the homological duals $a^*, b^*\in H_1(M/\partial M)$, two
disjoint, connected 
arcs beginning and ending in $\partial M$; slicing $M$
along these arcs gives a simply-connected polygon which can
serve as fundamental domain (for the action of $\pi_1(M)$ on the
universal cover); for our purposes we will label $a$ the edge
keeping the orientation of $a^*$, and $A$ its opposite edge with
the opposite orientation (see \figc{identifications});
similarly for $b$ and $B$.
Lifting a curve in $M$ to this fundamental polygon means representing
the curve as a set of arcs-with-identifications; each of these
{\em curve segments}  leads from
one of the edges $a, b, A, B$ to another; the orientation of
the curve defines a cyclic word in the four symbols: one records the
positive intersections as they occur. By construction, this word
represents the free homotopy class of the curve under consideration.

A curve segment is a {\em transversal} if it joins opposite   
edges of the fundamental domain, and  a {\em corner} otherwise.
Transversals correspond to consecutive $aa, AA, bb, BB$ in the
word; other combinations give corners.  Two corners are
{\em opposite} if they are diagonally opposed. Thus  $ab,ba$ and
$ab,AB$ correspond to diagonally opposed corners;
$ab$ and $ba$ have the \emph{same orientation}  whereas $ab$ and $AB$
have {\em reversed orientations}. In Figure \ref{identifications} the
curve $\langle baBBAba\rangle$ has two $ba$ corner segments diagonally
opposed to an $ab$ (same orientation) and a $BA$ (reversed orientation);
one $aB$ corner diagonally opposed to
opposed to an $ab$ (same orientation) and a $BA$ (reversed orientation);
one $aB$ corner diagonally opposed to
an $Ab$ corner (reversed orientations) and one $BB$ transversal.

A curve with only transversal self-intersections, and with
the smallest number of self-intersections
for its homotopy class (multiple points count with multiplicity:
a multiple intersection of $n$ small arcs counts as 
$\binom{n}{2}$ intersections)
is said to be {\em tight} (compare \cite{DT}: ``taut'').

Two-component multi-words $[w,w']$ enter into the surgery process.
We  define the {\em intersection number} $\inn (w, w')$ of two reduced 
cyclic words $w,w'$ to be the minimum number of intersections between 
a general-position  curve representing $w$ and  one representing $w'$.
The self-intersection number of the multi-word $[w,w']$ is then 
$\si ([w,w']) = \si (w) + \si (w') +\inn (w,w')$, and a pair 
of curves with that smallest number of self-intersections is also said 
to be {\em tight}. We also extend the $\alpha$ and $\beta$
notation to multiwords: $\alpha(]w',w''])$ is the total number
of occurrences of $a$ or $A$ in $w$ and $w'$; $\beta([w',w''])$ 
the total number of occurrences of $b$ or $B$.

\subsection{The surgery}

Whenever a cyclic word $w$ contains a pair of opposite corners,
it may be cut in two places, once in the middle of each of the
corners, to give two linear words.
These two linear words may be  reassembled (the corners
themselves are reassembled into transversals) into either
a new word $w'$ or a new multi-word $[w', w'']$
(according to the relative orientation of the corners);
if a multi-word $[v', v'']$ contains a pair of
opposite corners, one in each component, the two corners may be cut and
reassembled into two transversals yielding a new single word $v$.

\begin{figure}[httb]
\centerline{\includegraphics[width=7in]{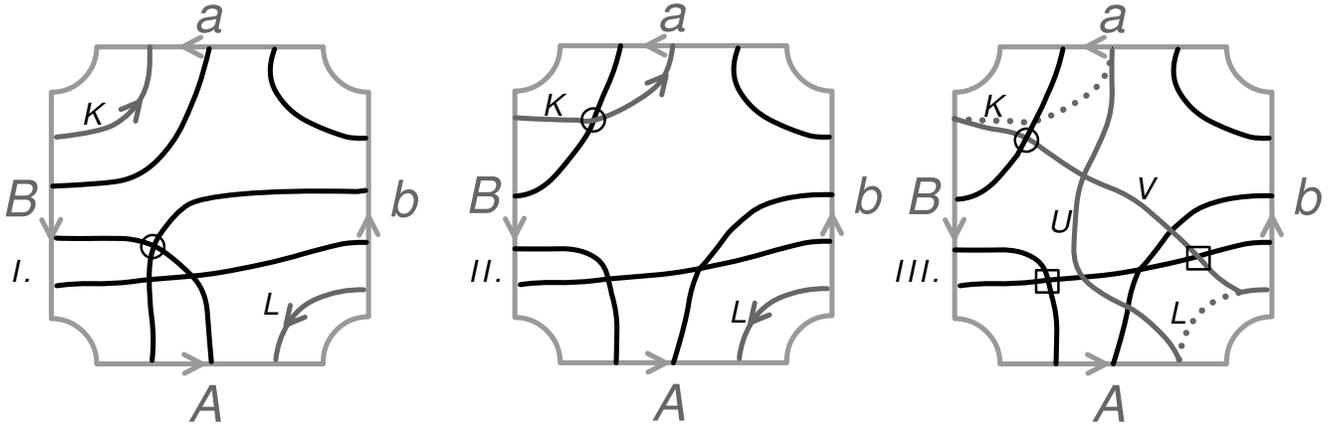}}
\caption{The cross-corner surgery $\langle baBB|Ab|a\rangle
\rightarrow \langle baBB|Ba|a\rangle$ as carried out on 
a tight representative curve. I. The corner $K$ corresponds to $b|a$;
$L$ corresponds to $B|A$. Note that the extension of $K$ intersects
that of one of its parallel corners (circled intersection). II. Before
surgery, that intersection is ``pushed,'' using a Rademeister-type-III
move, into the center of the surgery. III. $K$ and $L$ are excised,
$U$ and $V$ sewn in. The circled intersection
migrates to an intersection with $V$. The intersection
of $V$ with the original $BB$ spans a bigon with one of the original
vertices (squared intersections). $\si(\langle baBBBaa\rangle) = 6$.}\label{surg-example}\end{figure}

For a picture of the surgery on a curve, see \figc{surg-example};
in terms of the words,
the cutting and reassembly take one of the following forms:
\begin{align}
\langle \xx r| s\yy s| r\zz  \rangle & \rightarrow [\langle \xx r|
 r\zz \rangle,\langle   s\yy s| \rangle] \\
\langle  \xx r| s\yy R| S\zz  \rangle & \rightarrow \langle
\XX r| r\YY S| S\zz \rangle \\
[\langle \xx r| s\yy \rangle,\langle  \zz s| r\ww \rangle]
& \rightarrow \langle \xx r| r\ww\zz s| s\yy \rangle \\
\{\langle \xx r| s\yy \rangle, \langle \zz R| S\ww \rangle\}
& \rightarrow \langle \xx r| r\ZZ\WW s| s\yy \rangle
\end{align}
where 
${\mathsf x},
 \yy , {\mathsf z}, {\mathsf w}$ are
arbitrary (linear) subwords, and $R = r^{-1}, S = s^{-1}, 
{\mathsf X} = {\mathsf x}^{-1}$, etc. \
\start{defi}{Definition 2}
This cutting and reassembly
are called {\em cross-corner surgery} on the word $w$ or the
multi-word $[v', v'']$.
\end{defi}

It seems natural that transversals should contribute, more
than corners, to the self-intersection number of a curve.
Proposition \ref{cross corner}
makes this quantitative by showing that cross-corner
surgery, which eliminates two corners and adds two transversals,
always increases the self-intersection number by at least one. 
\start{prop}{cross corner}
 \begin{numlist}
\item[\rm{(1)}] If a word $w$ contains a pair of opposite corners with
reversed orientation
then cross-corner surgery will produce a new word $w'$, 
with $\alpha(w') = \alpha(w),~ \beta(w') = \beta(w)$, 
with one less block-pair, and with $\si(w')\geq \si(w)+1$.
\item[\rm{(2)}]  If a word $w$ contains a pair of opposite corners with
the same orientation  then cross-corner surgery will produce a multi-word 
$[w',w'']$, 
with $\alpha([w',w'']) = \alpha(w)$, $\beta([w',w'']) = \beta(w)$,
with one less block-pair, and with $\si([w',w''])\geq \si(w)+1$.
\item[\rm{(3)}]  
If a multi-word $[v',v'']$ contains a pair of opposite corners,
one in each component, irrespective of orientation, then cross-corner
surgery will produce a single word $v$  
with $\alpha(v) = \alpha([v',v'']),~ \beta(v) = \beta([v',v''])$,
with one less block-pair, and with $\si(v) \geq \si([v',v''])+1$. 
\end{numlist}
\end{prop}

This proposition is stated in terms of words, but its proof, given in the next subsection,
works by examining curves representing the words before
and after surgery; we first must fix a topological procedure
for carrying out cross-corner surgery on a curve. Specifically,
given a tight curve, or a tight pair of curves, representing
the candidates $w$ or $[v',v'']$ for cross-corner surgery,
we need to establish a systematic way of generating
curves representing the result $w'$, $[w',w'']$ or $v$ of
the surgery. We do this as follows:

\start{defi}{ccs-curves}{\em Cross-corner surgery on curves.} 
Suppose $r|s$ and $s|r$ or $R|S$ are the loci (that is, two diagonally 
opposite corners) in the word $w$ (or
the multi-word $[v', v'']$) chosen for surgery, and let $K$ and $L$
be the corresponding corners in a tight representative 
(see \figc{surg-example}). 
\begin{itemize} 
\item Preparation for the surgery.
If the extension of
any corner segment of the same type as $K$
(i.e. corresponding to the same letter sequence $rs$ or to the inverse
sequence $SR$) intersects the extension of $K$ in either direction 
before diverging, 
the curve is prepared for surgery  by a homotopy sliding that
(necessarily single) intersection onto the segment $K$ itself.
This deformation may be carried out by a sequence of Reidemeister-type-III
moves without changing the total number of intersections (see \figc{surg-example}, I and II). A similar
operation is carried out on the corner $L$.
\item Cutting and Sewing.
Corresponding to the
word permutation, corners $K$ and $L$ are removed and
replaced by transversals
$U$ and $V$. More
precisely, a line is drawn from a point on $K$ to a point on $L$, in
general position with respect to the rest of the curve, and cutting
any segment no more than once; that
line is  expanded into an $\mathcal{X}$-junction: $U$ routes the right
edge  of $K$ to the left edge of $L$, and vice-versa for $V$. 
\end{itemize}
\end{defi}

\subsection{Proof of \propc{cross corner}}

We will
obtain a lower bound on the increase
 in self-intersection number by counting the
vertices added and those possibly annihilated by the surgery. Annihilation
occurs through the creation of a {\em bigon:}  an 
immersed planar polygon  with
two vertices, and two edges with disjoint preimages  (a ``singular 2-gon''  in \cite{hs});  the bigon defines a homotopy of the
curve leading to the disappearance of its two vertices. An intersection
will be called {\em stable} if it is not the vertex  of a bigon; a curve
is tight if all its self-intersections are stable \cite{hs}.

\start{lem}{no-bigons} Cross-corner surgery does not
create any bigons spanned by a pair of pre-surgery vertices. 
\end{lem}

\begin{figure}[httb]
\centerline{\includegraphics[width=2in]{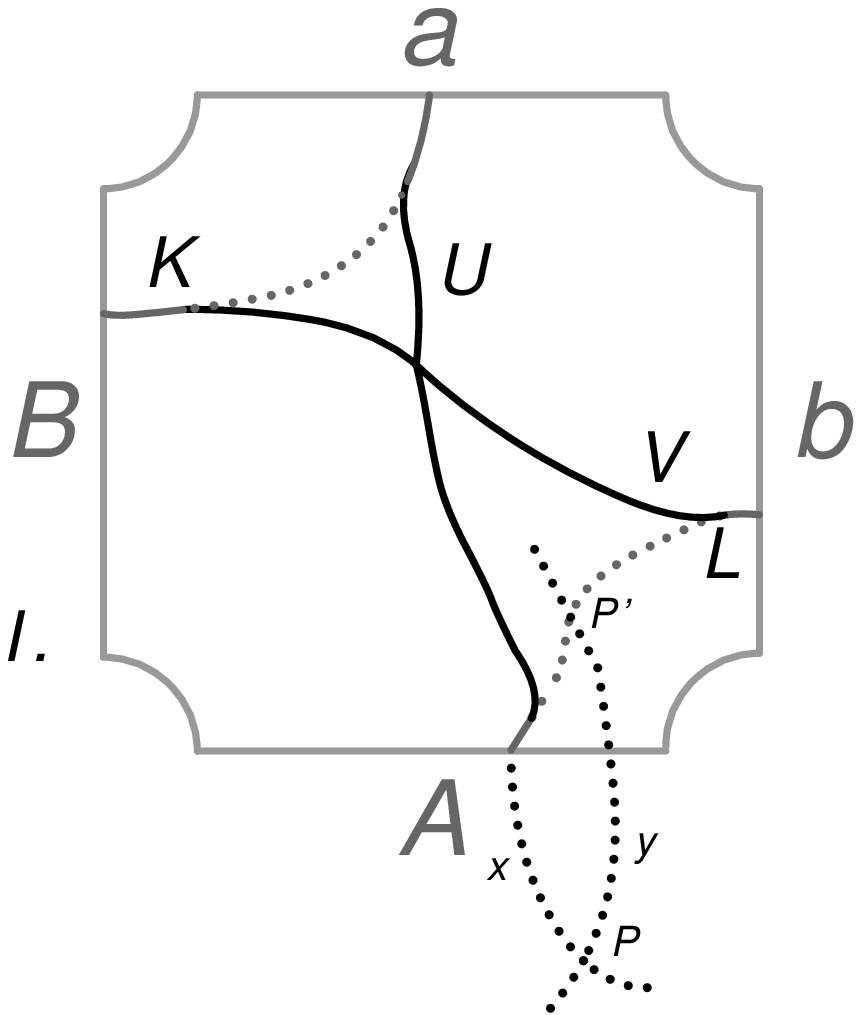}~~~~\includegraphics[width=2in]{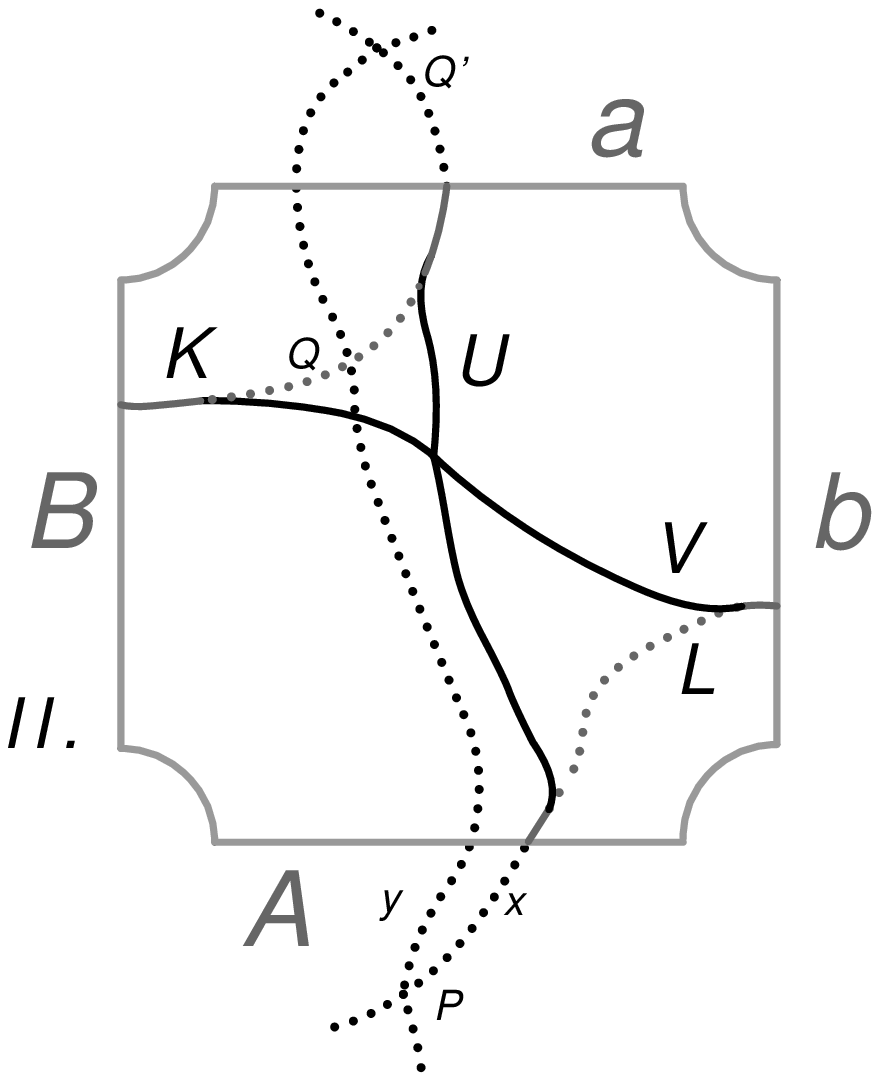}}
\caption{Cross-corner surgery cannot produce a bigon linking two pre-surgery vertices.}\label{n-o-b}
\end{figure}

\begin{proof} Since the initial curve is tight, the only way a pair $x, y$ of curve-portions starting from
a pre-surgery (``old'') vertex $P$ can lead to a bigon with another
old vertex is if one of those curve-portions  (say, $x$) 
contains one of the new segments $U$ or $V$, say $U$. Suppose the other one, i.e. $y$,  enters
inside the corner ($L$ in \figc{n-o-b}). Then (\figc{n-o-b}, I)
as $y$ follows
$x$ across the frame $y$ must intersect $L$ in an old vertex $P'$
cancelling $P$, contradicting tightness of the original curve. 
So $y$ must enter
outside $L$; then running parallel to $U$ across the frame 
 it must intersect the
opposite corner $K$ in an old vertex $Q$ (\figc{n-o-b}, II). 
Now if $x$ and $y$ meet in an old vertex $Q'$ so as to form a bigon cancelling $P$, then $Q'$ and $Q$
will span an old-vertex bigon. By tightness, this will require another use of the new segments.
Since each of $U$ and $V$ can only be used once by each of $x$ and $y$,
after at most four passes through the frame all the possibilities will be
exhausted; no such bigon can exist.

\end{proof}  

\begin{prooftext} {Proof of \propc{cross corner}} 
The curve surgery described in \defic{ccs-curves} yields one word if it is applied to a word that contains a pair of opposite corners with reversed orientation, a multi-word if it is applied to a word that contains a pair of opposite corners with the same orientation and a single word if it is applied to a multi-word that contains a pair of opposite corners, one in each component, irrespective of orientation. Thus, to prove (1), (2) and (3) it is enough to prove that in a  cross-corner surgery, the
number of new vertices minus the number of vertices cancelled by
new bigons is greater than or equal to one. We start by
classifying the new
vertices introduced by the surgery and the possible bigons in which they may
participate.

{\em Vertices:}

\begin{figure}[htp]
\centerline{\includegraphics[width=3in]{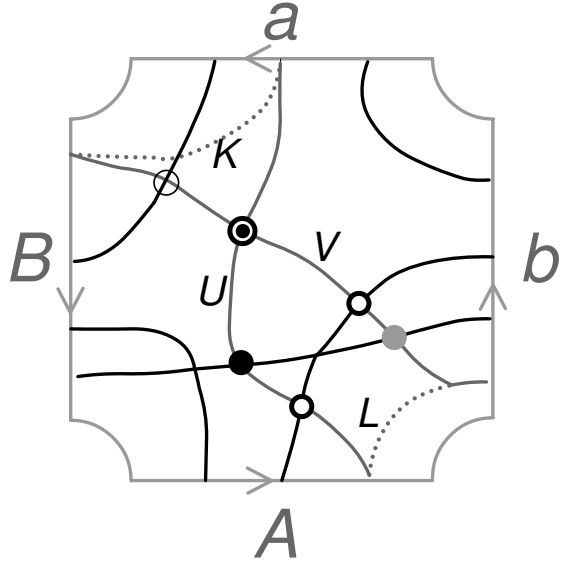}}
\caption{The new vertices created by a cross-corner surgery.}
\label{cross corner surgery}
\end{figure}

The surgery creates three
types of new vertices, shown as black, grey and white in
Figure \ref{cross corner surgery}, as follows.
\begin{numlist}
\item[{\rm (1)}] (black) Stable intersections between $U$ and horizontal transversals (i.e
segments corresponding to $bb$ or $BB$ in the initial word $w$), between
$V$ and vertical transversals, and (bullseye)
the stable intersection between $U$
and $V$. 
 
\item[{\rm (2)}] (grey) Intersections between $U$ and other vertical transversals, and
between $V$ and other horizontal transversals. These are potentially
vertices of bigons. 

\item[{\rm (3)}] (white) Intersections between $U$, $V$ and remaining corner segments. In 
\figc{cross corner surgery}
only those of type $ab$ or $BA$ are shown; there is typically another
family in the opposite corner corresponding to types $ba$ or $AB$. These
are also potentially vertices of bigons.

\item[{\rm (4)}] In addition, 
the circled vertices in Figure \ref{cross corner surgery} 
are those inherited by the new curve from the old.  
These correspond to the intersections between $K$  or $L$ with
other corners of the same type; such a corner is labeled $J$
in \figc{yy}. 
Focussing on $K$,
let us label $x$ and $y$ the two ends of the segment $K$, and
by $w$ the intersection point of the new segments $U$ and $V$.
The segment $K$ and the broken curve $uwv$ are fixed-endpoint homotopic;
it follows that for any original segment having exactly one endpoint between
$u$ and $v$ (e.g. $J$)   
 the intersection with $K$ will migrate to an intersection with $U$ or $V$ during that homotopy (with $V$ if the outside end of $J$
is on the $B$ side --as in \figc{yy}-- and with $U$ if it is on the
$a$ side).
\end{numlist} 
\begin{figure}[httb]
\centerline{\includegraphics[width=3in]{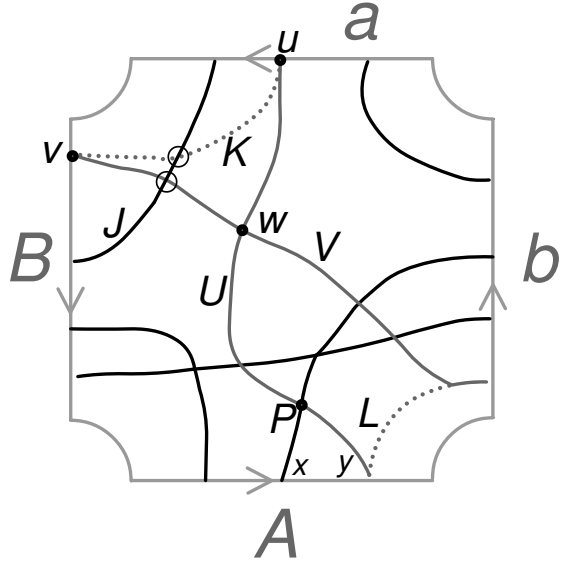}}
\caption{Vertices inherited by new curve from old; $P$ is an example of a type-3 vertex.}
\label{yy}
\end{figure}

{\em Bigons:}

The only bigons that need to be examined are those where
one of the spanning vertices is an old vertex or a type-4 vertex; 
because if two new vertices 
form a bigon
and cancel, that does not affect the inequality we need to prove. So,
letting 1, 2, 3, and 4  
represent vertices so labeled above,  letting  $x$, $x'$ represent
self-intersections of the original curve, and keeping in mind that type-1
vertices are stable, and that bigons of type $(x,x')$ cannot occur 
(\lemc{no-bigons}),
we need only examine bigons of type $(2,x)$, $(2,4)$,  $(3,x)$, $(3,4)$,  
$(4,x)$,  and $(4,4)$.

\begin{romlist}
\item
$(4,x)$ and $(4,4)$. A vertex of type 4 can
span a bigon in only  one of its quadrants; but in that
quadrant a bigon would imply a bigon with the old vertex
from which the type-4 vertex was inherited;
so a $(4,x)$ would imply an $(x',x)$, and a $(4,4)$ would
imply a $(4,x)$; so neither $(4,4)$ nor $(4,x)$ can
occur.
\item
 $(2,x)$ and $(2,4)$.
A type-2 vertex $y$ may span a bigon with an old vertex $x$; the type-2
vertex is either the intersection of $U$ with another vertical
transversal, or $V$ with another horizontal. In the first case
(the second case is similar),
that vertical transversal must also intersect $V$, creating a new
(type-1) stable intersection $z$.  In total we will have added two vertices
($y$ and $z$), and lost two vertices ($y$ and $x$) to a bigon.
The inequality is not affected. 

Since, arguing as in (i), a $(2,4)$ bigon would imply a $(2,x)$
bigon, the loss of the 4 would be balanced by the gain of the
corresponding new type-1 vertex, and again the inequality would not be
affected. 

\item $(3,x)$ and $(3,4)$. Figure \ref{yy} shows a typical type-3 vertex $P$. It can only span a bigon in one quadrant; label $x$
and $y$ the two segments issuing from $P$ in that direction. Because of the way the curve is prepared for surgery, $x$ and $y$ cannot
be continued with old segments to form a bigon cancelling $P$.
We need to discuss the
possibility that after surgery their extensions could incorporate
$U$ or $V$ or both and then form such a bigon. This $P$,
$x$ and $y$ exactly match the notation of \lemc{no-bigons}; and
the proof of that lemma applies here as well: no such bigon can exist. Since a $(3,4)$ bigon would imply a $(3,x)$ bigon,
no type-$(3,4)$ bigons can exist either.

\end{romlist}

In summary, cross-corner surgery generates one special stable vertex
(the intersection of $U$ and $V$) plus other new vertices
of types 1, 2, 3, and displaced vertices of type 4. 
Vertices of type 1 are stable. Some of the vertices of type 2 and 3
form bigons with each other 
and cancel out. Vertices of type
 3 and 4 cannot form bigons with
pre-surgery vertices, and any old or type-4 vertex cancelled by a type-2 vertex
can be replaced in the count   by the corresponding type-1 vertex.
It follows that cross-corner surgery increases the self intersection 
number by at least one.

\end{prooftext}

\section{Linked pairs}\label{lpsect}

Ultimately the calculation of $\si(w)$ or $\inn([w',w''])$
 can be made directly from $w$ or $[w',w'']$, by
counting {\em linked pairs}. In this section we give a
simplified definition appropriate for the punctured torus,
we list two theorems from \cite{chas} giving the correspondence
between linked pairs and intersection points, and we summarize
explicit calculations of intersection and self-intersection
numbers for certain families of words with a small number of
block-pairs. Linked pairs have also been defined and studied
by Cohen and Lustig \cite{cl}.

\start{nota}{sym} From now on,  we will
use the symbols $p,  q,  r$,  $s,  p_1,  q_1$, etc. to represent
letters from the alphabet $a, b, A, B$, with $P = p^{-1}$,
etc. The symbols $v, w, v', w'$, etc will represent cyclic
words in that alphabet, e.g. $w = \langle abbaB\rangle  = \langle aBabb\rangle$. 
Sans-serif
symbols $\uu, \vv, \yy$ will represent linear words in the
alphabet $\{a,b,A,B\}$ with $ \rr, {\mathsf s}, \RR, {\mathsf S}$ representing
homogeneous blocks of letters  $rr\dots r, ss\dots s, RR\dots R,
SS\dots S$ respectively.
As before, $\VV = \vv^{-1}, \RR = \rr^{-1}$, etc. 
\end{nota}

\begin{proclama}{Orientation:} \emph{For these purposes 
we identify the boundary
of our fundamental domain with a clock face, with $a,b,A,B$ at
3, 6, 9 and 12 o'clock. Given six letters $p, q, r, p', q', r'$
from the alphabet $a, b, A, B$, we say that the triples $p, q, r$
and $p', q', r'$ are {\em similarly oriented} if the arcs 
$pqr$ and $p'q'r'$ have the same orientation on the clock face. This
implies that the three points in each triple are distinct. }
\end{proclama}

\start{defi}{def} Let $w$ (resp. $[w',w'']$) be a 
primitive reduced cyclic word (resp. a multi-word with primitive reduced 
cyclic components), corresponding to a
free homotopy class (resp. a pair of free homotopy classes) of curves 
on the punctured torus.  

\begin{figure}[httb]
\centerline{\includegraphics[width=3in]{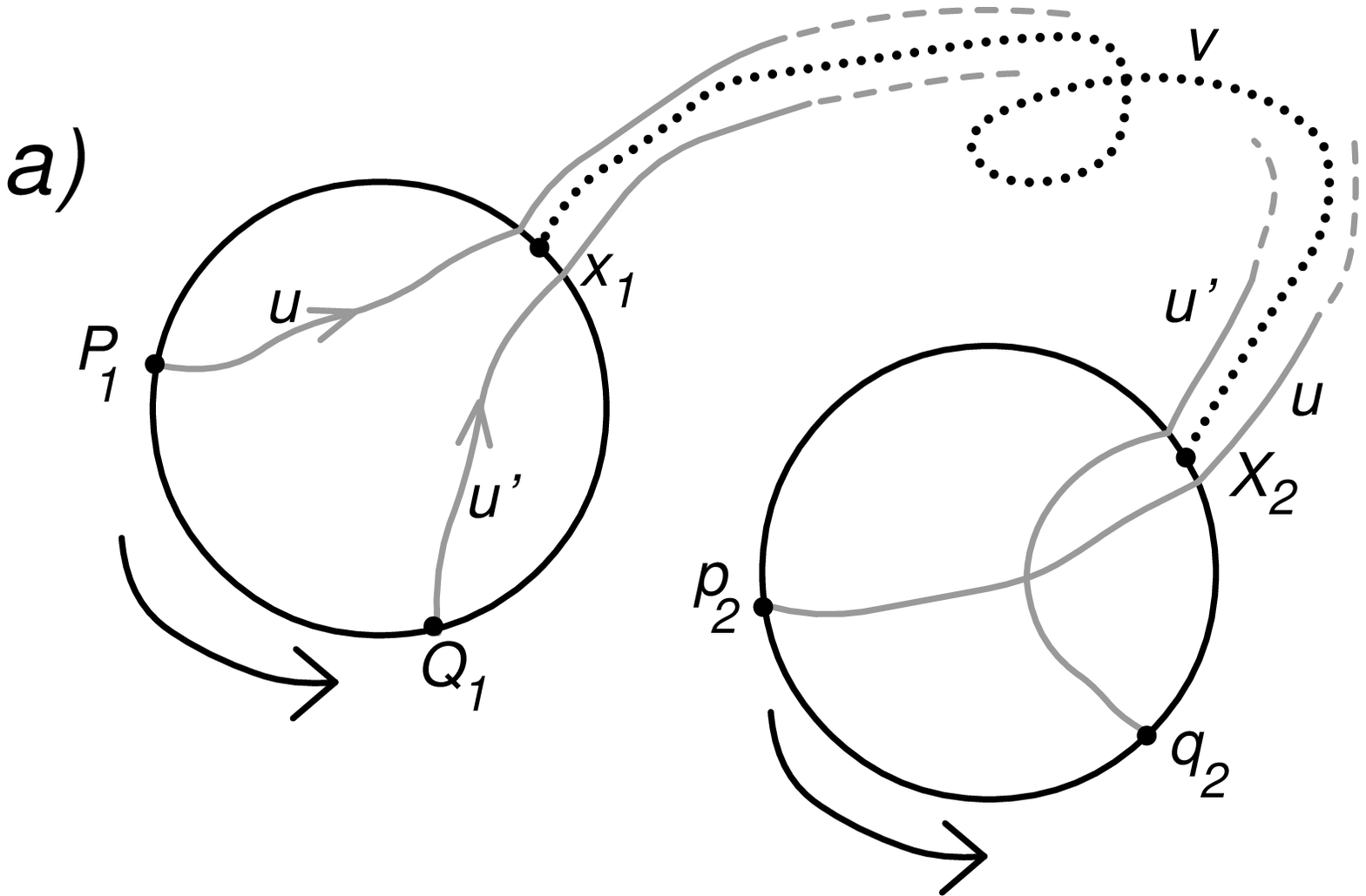}~~\includegraphics[width=3in]{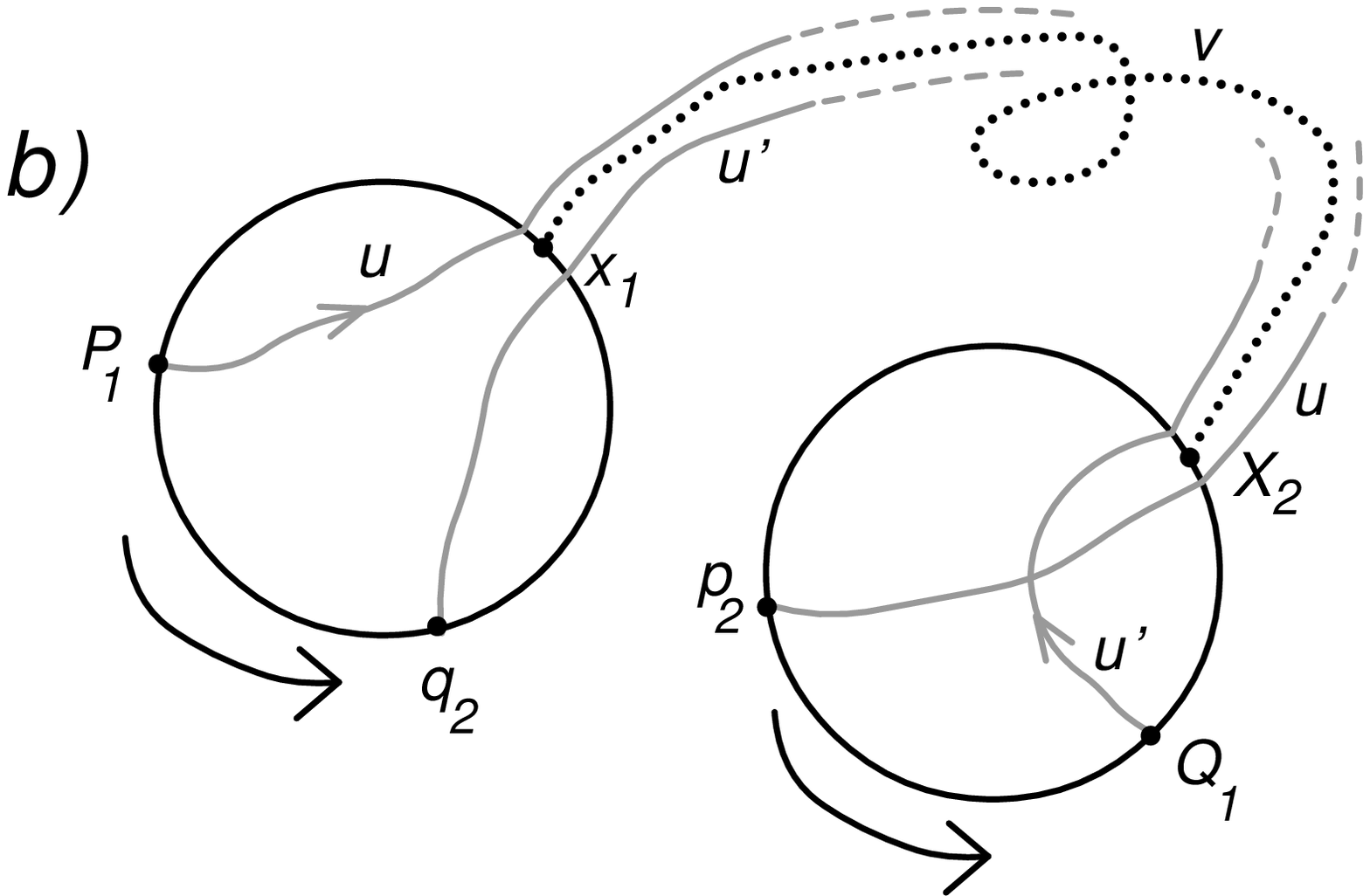}}
\caption{Linked pairs. Here $ \yy  = x_1{\mathsf v}x_2$. 
a) The linked pair
is ($p_1 \yy p_2$, $q_1 \yy q_2$). Since the orientations
of $(P_1,Q_1, x_1)$ and $(p_2,q_2, X_2)$ are the same, the
curve segments must intersect.  $~~$ b) The linked pair
is ($p_1 \yy p_2$, $q_2{\mathsf Y}q_1$). Since the orientations
of $(P_1,q_2, x_1)$ and $(p_2,Q_1, X_2)$ are the same, the
curve segments must intersect.} \label{pair-fig} 

\end{figure}

Let $\uu'$ and $\uu''$ be two possibly
overlapping but distinct {\em linear}
subwords, both of the same length $\geq 2$, of $w$ or of $[w',w'']$ 
(in that case let $\uu'\subset w'$, and let $\uu''\subset w''$).
The pair of words $\{\uu',\uu''\}$ is a {\em linked
pair} if one of the following criteria is satisfied (see  \figc{pair-fig}).
\begin{lplist}
\item $[\uu',\uu'']$ is one of the following words: $\{ aa, bb\}, \{ aa, BB\},\{  AA, bb\}, \{AA, BB\}$. 
\item  \begin{romlist}
\item (length 3) $\uu' = {\it p}_1{\it r}{\it p}_2$, $\uu'' = 
{\it q}_1{\it r}{\it q}_2$ (same $r$) with
$P_1Q_1r$ and $p_2q_2R$ similarly oriented, or 
\item (length $n$) $\uu' = {\it p}_1\yy {\it p}_2$, $\uu'' = 
{\it q}_1\yy {\it q}_2$, 
$\yy = {\it x}_1\vv {\it x}_2$ ($\vv$ possibly empty) with
$P_1Q_1x_1$ and $p_2q_2X_2$ similarly oriented.
\end{romlist}
\item \begin{romlist}
\item (length 3) $\uu' = {\it p}_1{\it r}{\it p}_2$, $\uu'' = 
{\it q}_1{\it R}{\it q}_2$ ($R=r^{-1}$) with
$P_1q_2r$ and $p_2Q_2R$ similarly oriented, or
\item (length $n$)
$\uu' = {\it p}_1\yy {\it p}_2$, $\uu'' = 
{\it q}_1\YY {\it q}_2$, 
$\yy = {\it x}_1\vv {\it x}_2$ ($\vv$ possibly empty) with
$P_1q_2x_1$ and $p_2Q_1X_2$ similarly oriented.
\end{romlist}
\end{lplist}
\end{defi}

\start{rem}{xxx}  $\{{\mathsf u},{\mathsf u'}\}$ is a linked pair of type
(2) if and only if
 $({\mathsf u},{\mathsf U'})$ is a linked pair of type (3).
\end{rem}

Tables~\ref{tab-a}, \ref{tab-ab} and \ref{tab-aba} summarize for future reference the pairing between
various subwords of a cyclic word $w$. In these tables an ``='' means
that the row word and the column word have the same first or last letter 
(so
they cannot form a linked pair); ``N'' means that there is no end
matching but that the pair fails the orientation criterion; ``Y''
means that the row word and the column word form a linked pair.

\begin{table}[htbp]
   \centering
  \begin{tabular}{|l|ccccccc| }
\hline
    Words    & $a^{i+1}b$ & $ba^{i+1}$ & $a^{i+1}B$ &$Ba^{i+1}$&$ba^ib$&$a^{i\
+2}$&$Ba^{i}B$ \\\hline

    $a^{i+1}b$ & = &N&= & Y&=&=&Y\\
     $ba^{i+1}$ &N &=&Y&=&=&=&Y\\
     $a^{i+1} B$ &=&Y&=&N&Y&=&=\\
     $Ba^{i+1}$&Y&=&N&=&Y&=&=\\
     $ba^ib$&=&=&Y&Y&=&Y&Y\\
     $a^{i+2}$&=&=&=&=&Y&=&Y\\
     $Ba^{i}B$&Y&Y&=&=&Y&Y&=
 \\\hline     
   \end{tabular}
   \caption{Linking of pairs of words with $\YY={\it a^{i}}$ 
(notation from \defic{def}, =, Y, N explained in the text).} 
 \label{tab-a}
\end{table}

\begin{table}[htbp]
   \centering
  \begin{tabular}{|c|cccc| }
\hline
      Words    & $aa^ib^{j}b$ & $ba^{i}b^{j}a$ & $ba^ib^{j}b$ & $aa^{i}b^{j}a$ \
\\\hline
    $aa^ib^{j}b$&=&N&=&=\\
 $ba^{i}b^{j}a$&N&=&=&=\\
   $ba^ib^{j}b$&=&=&=&Y\\
     $aa^{i}b^{j}a$&=&=&Y&=

    \\\hline
   \end{tabular}
   \caption{Linking of pairs of words with $\YY={\it a^{i}b^{j}}$  
(notation as  in \defic{def}).}\label{tab-ab}

\end{table}

\begin{table}[htbp]
   \centering
  \begin{tabular}{|c|cccc| }
\hline
 Words    & $aa^ib^{j}a^{k}b$ & $ba^{i}b^{j}a^{k}a$ & $ba^ib^{j}a^{k}b$ & 
$aa^{i}b^{j}a^{k}a$ \\\hline
    $aa^ib^{j}a^{k}b$&=&N&=&=\\
$ba^{i}b^{j}a^{k}a$&N&=&=&=\\
   $ba^ib^{j}a^{k}b$&=&=&=&Y\\
     $aa^{i}b^{j}a^{k}a$&=&=&Y&=
    \\\hline
   \end{tabular}
   \caption{Linking of pairs of words with $\YY={\it a^{i}b^{j}a^k}$  
(notation as  in \defic{def}).}\label{tab-aba}
\end{table}

The following theorem will be used  to compute the  self-intersection numbers of certain words and multi-words (see \propc{lp} and Appendix~\ref{appendix}). This theorem is a direct consequence of \cite[Theorems 3.9 and 3.10 and Remarks 3.10 and 3.11]{chas} and \cite[Theorem 3.12 and Remark 3.13]{chas}.

 \start{theo}{clb si}  
Let $v$ and $w$ be a cyclic reduced words in the alphabet $\{a,b,A,B\}$. 
Suppose that $w=\langle u^k\rangle$ is the $k$th power ($k\geq 0$)
of the primitive reduced cyclic word $u$.
\begin{numlist}
\item[{\rm (1)}]  If $k=0$, so $w$ is primitive, $\si(w)$ is equal to  
the {\em number of linked pairs of $w$}, i.e. the cardinality of 
the set of unordered pairs $\{{\mathsf u}, {\mathsf u'}\}$, 
${\mathsf u}$ and ${\mathsf u'}$ linear subwords of $w$,
with ${\mathsf u}$ and ${\mathsf u'}$ linked as in \defic{def}.  

\item[{\rm (2)}] In general, $\si(w)$ is less than or equal to  $(k-1)$ plus  
the number of linked pairs of $w$.

\item[{\rm (3)}] $\inn(\{v,w\})$ equals the number of ordered pairs $({\mathsf u},{\mathsf u'})$ for which there exist positive integers $j$ and $k$ such that ${\mathsf u}$  is an occurrence of a subword of $v^j$, but not a subword of $v^{j-1}$, ${\mathsf u'}$ is an  occurrence of a subword of  $w^k$, but not a subword of $w^{k-1}$, and ${\mathsf u}$, ${\mathsf u'}$ are linked as in \defic{def}. (See \remc{k}.)

\end{numlist}
\end{theo}
In this work, only the following simple instances of \theoc{clb si}(3) will be necessary.

\start{rem}{k} \begin{numlist}\item[{\rm (1)}] $\inn( \langle a^ib^j \rangle,  \langle a^k\B^l \rangle)$ equals the number of ordered pairs $({\mathsf u},{\mathsf u'})$  such that ${\mathsf u}$  is an occurrence of a subword of $\langle a^ib^j \rangle$,  ${\mathsf u'}$ is an  occurrence of a subword of  $\langle a^k\B^l \rangle$ and ${\mathsf u}$, ${\mathsf u'}$ are linked as in \defic{def}.
\item[{\rm (2)}] $\inn( \langle a^ib^ja^kb^l\rangle, \langle a^m\B^n \rangle)$ equals the number of ordered pairs $({\mathsf u},{\mathsf u'})$  such that ${\mathsf u}$  is an occurrence of a subword of $\langle a^ib^ja^kb^l\rangle$,  ${\mathsf u'}$ is an  occurrence of a subword of  $ \langle a^m\B^n \rangle$ and ${\mathsf u}$, ${\mathsf u'}$ are linked as in \defic{def}.
\end{numlist}
This is because if  $[v,w]=[\langle a^ib^j \rangle,  \langle a^k\B^l \rangle]$ or  $(v,w)=[\langle a^ib^ja^kb^l\rangle, \langle a^m\B^n \rangle]$, $J$ and $K$ are non-negative integers, $\mathsf{u}$ is a linear word which is an occurrence of a subword of $v^{J}$ and $w^{K}$ then $\mathsf{u}$  is an occurence of a subword of $v$ and $w$.  
\end{rem}

In principle the self-intersection number corresponding to
any particular word can be ascertained combinatorially by a
count of linked pairs. The number of steps in this calculation
increases rapidly with the length of the word, but
it can be carried out completely for words with a small number of
block pairs. The results of these calculations are given in
\propc{lp}, with the work itself presented in Appendix A.

\start{prop}{lp}

\begin{numlist}
\item[{\rm (1)}] $\si( \langle a^ib^j \rangle)=(i-1)(j-1)$.
\item[{\rm(2)}]
\begin{equation}\label{siabab}
\si(\langle a^ib^ja^kb^l\rangle)  \left \{ \begin{array}{ll}  
\le (i +k - 2)(j + l -2) + 1 & \mbox{if  $k = i$ and $l = j$,} \\
= (i + k - 2)(j + l - 2) + |i -k| + |j - l| - 1 &   \mbox{otherwise.}
\end{array} \right . 
\end{equation} 
\item[{\rm(3)}] $\si( \langle a^ib^ja^k\B^l \rangle)=(i+k-1)(j+l-1)$.
\item[{\rm(4)}] $\si(\langle a^{i}b^{j}A^{k}B^{l} \rangle=(i+k-1)(j+l-1)-1$.
\item[{\rm(5)}] $\si( \langle a^ib^ja^kb^la^mB^n \rangle)=(i+k+m-1)(j+l+n-1) -2    (k+\min(j,l)-1)$.
\item[{\rm(6)}] $\inn( \langle a^ib^j \rangle,  \langle a^k\B^l \rangle)=il+kj$.
\item[{\rm(7)}] $\inn( \langle a^ib^ja^kb^l\rangle, \langle a^m\B^n \rangle)=(i+k)n+m(j+l)$.
\end{numlist}
\end{prop}

\start{cory}{lpcory}
$$\si(\langle a^ib^ja^kb^l\rangle)  \left \{ \begin{array}{ll}  
= 1 & \mbox{if $i = k$ and $j = l$}\\
& \mbox{and  $i = 1$ or $j=1$}\\
\leq (i +k - 1)(j + l -1) - 4 & \mbox{if  $k = i \ge 2 $ and $l = j\ge 2$,} \\
\leq (i +k - 1)(j + l -1) -2 &   \mbox{if $i\neq k$ or $j\neq l$.}
\end{array} \right . $$
\end{cory}

\begin{proof}
 It follows from \propc{lp} (2) that
if $i=k$ and $j=l$ and either pair is 1, then the $\si \le 1$; and if
both are $\geq 2$ then 
\begin{align*}(i +k - 2)(j + l -2) + 1 &=
(i+k-1)(j+l-1)-(i+k-1)-(j+l-1)+2 \\
 &\leq (i+k-1)(j+l-1)-4.\end{align*}
If $i\neq k$ or $j\neq l$ then \propc{lp} (2)
gives 
\begin{align*}
\si( \langle a^ib^ja^kb^l \rangle)&=(i+k-2)(j+l-2)+|i-k|+|j-l|-1\\
 &=(i+k-1)(j+l-1) -(i+k) +1 -(j+l) +1 +1 +|i-k|+|j-l|-1\\
 &=(i+k-1)(j+l-1)-2\min(i,k)-2\min(j,l) +2\\
 &\leq (i+k-1)(j+l-1)-2.\end{align*}

\end{proof}

The next remark is useful in the proof of \propc{ababaB}. 

\start{rem}{corners} In the punctured torus, it follows from \defic{def} that if ${\mathsf P}
 = rs {\mathsf u}sR$, where $r$ and $s$ are distinct letters 
and ${\mathsf u}$ is an arbitrary linear word, then
$\{{\mathsf P},{\mathsf Q}\}$ is not a linked pair for any ${\mathsf Q}$.
\end{rem}

\section{Proof of \theoc{main-th}}\label{proof-ub}

\subsection{Detailed Strategy of Proof}

This subsection amplifies the sketch presented in Subsection \ref{sketch},
continuing with the notation from Definition \ref{def-block-pair}
and Section \ref{sec-[httb]}.

Given an arbitrary
reduced cyclic word, 
we prove that its self-intersection
number must be less than or equal to that
of a word of the same length with few enough block-pairs to be
 amenable to a linked-pair self-intersection-number
calculation. 

This ``amalgamate and conquer'' strategy is implemented by
 cross-corner surgery, which reduces the number
of block-pairs in $w$ while conserving $\alpha(w)$ and $\beta(w)$
 and increasing  
$\si(w)$.

\begin{figure}[httb]
\centerline{\includegraphics[width=5in]{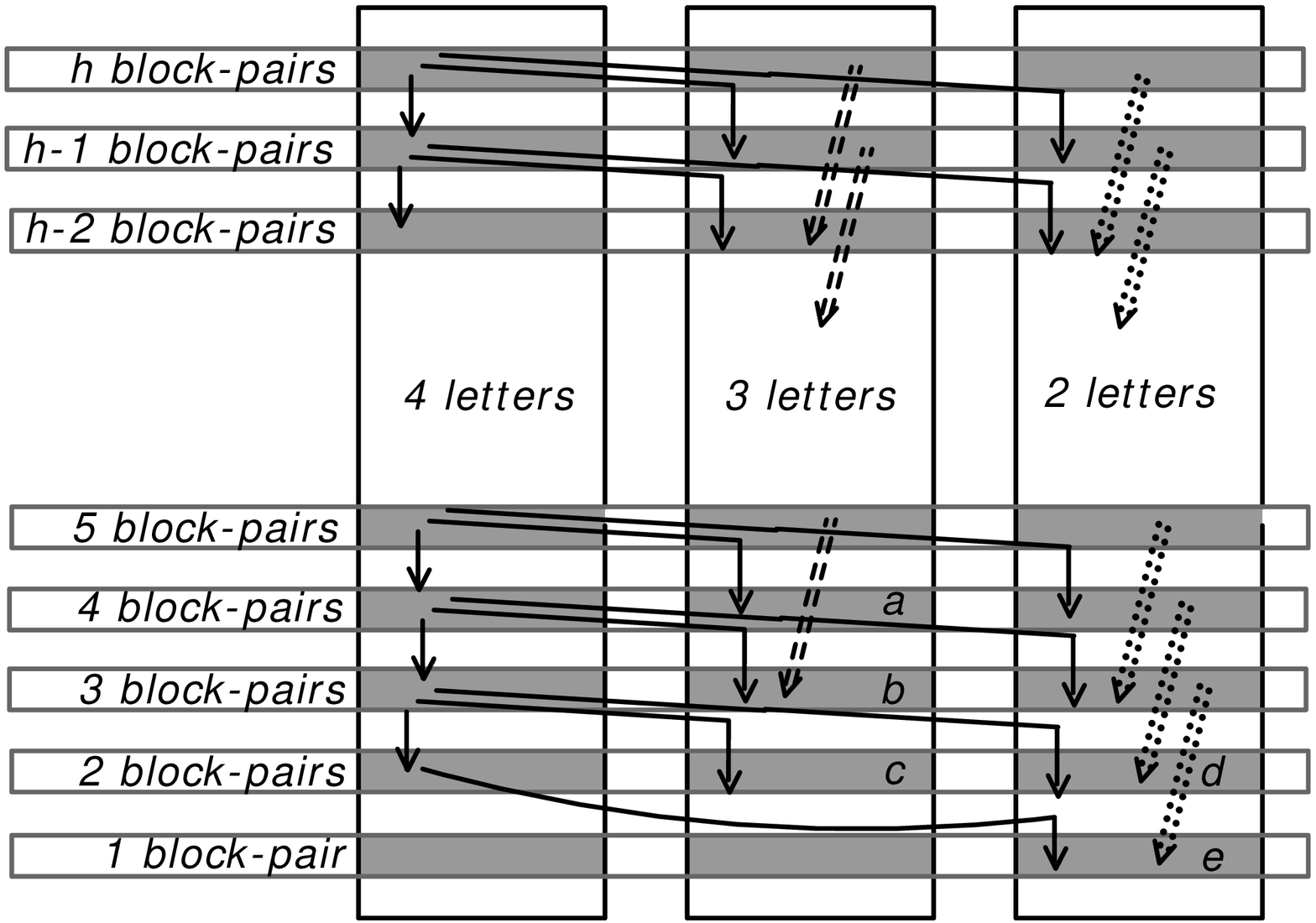}}

\caption{Flow chart of proof. The arrows correspond to
possible cross-corner surgeries. 
Straight-line-arrows: \lemc{four}; 
dashed arrows: \lemc{Three-letter 5b}; 
dotted arrows:  \lemc{two}. 
Terminal cases: (a) = \lemc{abababaB};
(b) = \propc{lp} (5);
(c) = \propc{lp} (3);
(d) = \coryc{lpcory};
(e) = \propc{lp} (1). 
Note that the self-intersection number
of a pure power (word with 1 letter)
can be calculated directly ($\si(r^k)=k-1$)
and then compared with the maximum for general words of
the same length; see \remc{rem1}.} \label{flow}
\end{figure}

The detailed procedure at each step in the reduction
depends on the number of different letters
occurring in the word  \figc{flow}). As we will see, 
\begin{itemize}
\item a word that uses
all four letters 
is always a candidate for cross-corner surgery using opposite
corners with reversed orientation; the result
will be a single word with one less block-pair
(since this surgery reverses the orientation of part of the word,
the number of different letters may change);
\item if a word uses
exactly three of the four letters and has at least five block-pairs,
or if it uses only two of the four letters and
has at least three block-pairs,
then two cross-corner surgeries will reduce the number of
block pairs by two (the intermediate stage is a two-component
multi-word) and increase the self-intersection number by at least two; 
these surgeries permute the letters in the word,
and so the new word still uses three letters or two letters if
the old one did. 
\end{itemize}
So the
words remaining are:
\begin{itemize}
\item words with three letters and 

\begin{alphlist}
\item 4 block-pairs
($\langle \mathsf{rsrsrsrS} \rangle $, 
$\langle\mathsf{rsrsrSrS}  \rangle $ and
$\langle \mathsf{rsrSrsrS}  \rangle $),
\item 3 block-pairs
($\langle \mathsf{rsrsrS} \rangle $), or  
\item 2 block-pairs 
($\langle \mathsf{rsrS}  \rangle $); 
\end{alphlist}

\item words with two letters and 

\begin{alphlist}
\item[(d)] 2 block-pairs ($\langle \mathsf{rsrs} \rangle $) or 
\item[(e)] 1 block-pair ($\langle \mathsf{rs} \rangle $). 
\end{alphlist}

\item pure powers. 
\end{itemize}


\subsection{Preparatory Lemmas}
In these lemmas and their proofs,  Notation \ref{sym} will be used.

\start{lem}{four} If a reduced cyclic word $w$
contains all four letters $a, \A, b, \B$ then 
there exists a word $w'$ with the same $\alpha$ and $\beta$ values,
with one less block-pair, and with $\si(w')\geq\si(w)+1.$
\end{lem}

\begin{proof} Claim: such a word must contain two corners with reverse orientation. In fact,   let $w$ be a reduced cyclic word which contains
all four letters (such a word must have at least two block-pairs) and
 which does not contain two subwords of the form $xy$ and $XY$,  
where  $x \in \{a,\A\}$ and  $y \in \{b,\B\}$ or vice-versa.  Now $w$ must contain at least one of $ab$ and $aB$; suppose $w$ contains $ab$.
Then $w$ does not contain $AB$. 
So every  $\BB$-block must be preceded by an  $a$. Since there is at
least one such block, $w$ must contain $a\B$, which implies that 
$w$ does not contain $\A b$.
Since $w$ does not contain $A\B$ or $\A b$, there is no letter possible 
after  an $\AfA$-block. Since there is at least one such block, our
hypothesis leads to a contradiction.

The lemma now follows from  \propc{cross corner},(1).
\end{proof}

\start{lem}{Three-letter 5b} Suppose a cyclic word $w$ uses
exactly three distinct letters from the set $\{a, \A, b, \B\}$
and has five or more block-pairs. Then there exists a word $w'$, with two fewer block-pairs, with the same $\alpha$ and
$\beta$ values, and such that  $SI(w')\geq \si(w)+2$.
\end{lem}
\begin{proof} 

Suppose the three letters are $a$, $b$ and $B$. The block-pairs are either 
$\afa\bb$'s or $\afa\BB$'s. We may suppose there are at least three 
$\afa\bb$'s. Hence $w$ has the form  
$\langle \afa\bb~\afa\uu~\afa\bb~\afa\vv~\afa\bb~\afa\yy \rangle $,
where $\uu,\vv,\yy$ represent (possibly empty) 
blocks of letters.

We pick two {\em consecutive} $\afa\bb$ block-pairs and apply 
\propc{cross corner} (Cross Corner Surgery) as follows:
$\langle \afa|\bb~\afa\uu~\afa\bb|\afa\vv~\afa\bb~\afa\yy \rangle $
$\rightarrow [\langle \bb~\afa\uu~\afa\bb|\rangle, 
\langle\afa|\afa\vv~\afa\bb~\afa\yy \rangle]$ 
$ =[\langle \bb\afa\uu\afa \rangle,
\langle \afa\vv\afa\bb\afa\yy \rangle]$ $= [v',v'']$.

We have lost one $ab$ corner and one $ba$ corner, so the number of block-pairs
has gone down by one. On the other hand, \propc{cross corner} 
guarantees that $\si([v',v'']) \ge \si(w)+1$.

Our consecutive corner condition guarantees that both $v'$ and $v''$
contain both $ab$ and $ba$, so the multi-word $[v',v'']$ is a candidate for
a second surgery, for example:
$[\langle |\bb\afa\uu\afa \rangle,
\langle \afa\vv\afa\bb|\afa\yy \rangle]$
$\rightarrow$ $\langle |\bb\afa\uu\afa|\afa\yy\afa\vv\afa\bb \rangle
= \langle \bb\afa\uu\afa\yy\afa\vv\afa \rangle = w'.$

We have lost another pair of corners, so the number of block-pairs has 
gone down by one more; \propc{cross corner} 
guarantees that  $\si(w') \ge \si([v',v''])+1$, and thus
$\si(w')\ge \si(w) +2$. The $\alpha$ and $\beta$ values are clearly the  same. 
\end{proof}

\start{lem}{abababaB} If $w$ has one of the forms
$\langle \afa\bb\afa\bb\afa\bb\afa\BB \rangle$,
$ \langle \afa\bb\afa\BB\afa\bb\afa\BB \rangle$,
$\langle \afa\bb\afa\bb\afa\BB\afa\BB \rangle$ 
then there exists a word 
$w'$ with two block-pairs, the 
same $\alpha$ and $\beta$ as $w$ and such that $\si(w') \ge \si(w)+2$.
\end{lem}
\begin{proof} 
\begin{enumerate}
\item $\langle \afa\bb\afa\bb\afa\bb\afa\BB \rangle$.
We apply 
cross corner surgery (\propc{cross corner}) as follows:
$\langle \afa\bb|\afa\bb\afa|\bb\afa\BB \rangle$
$\rightarrow$ $[\langle \afa\bb|\bb\afa\BB \rangle,
\langle |\afa\bb\afa ]$ $ = [\langle \afa\bb\afa\BB \rangle,
\langle \bb\afa ] = [v',v''].$ 
The multiword $[v',v'']$ has the same $\alpha$ and $\beta$ values as $w$. 
Furthermore, $\si([v',v'']) \geq \si(w) + 1$. 
Another application of \propc{cross corner}:
$ [\langle \afa|\bb\afa\BB \rangle, \langle \bb|\afa\rangle ]$
$\rightarrow$ $\langle \afa|\afa\bb|\bb\afa\BB \rangle = 
\langle \afa\bb\afa\BB \rangle =w'$
yields a word $w'$
 with two block-pairs, the same $\alpha$ and $\beta$ values as  $w$, 
and  $\si(w')\geq \si(w)+2$.

\item $ \langle \afa\bb\afa\BB\afa\bb\afa\BB \rangle$.
Cross-corner surgery
$ \langle \afa\bb\afa\BB|\afa\bb\afa|\BB \rangle$
$\rightarrow$ $[\langle \afa\bb\afa\BB|\BB \rangle,
\langle \afa\bb\afa|\rangle]$ $=[\langle \afa\bb\afa\BB \rangle,
\langle \bb\afa\rangle]$ leads to the same half-way step as the
previous case.

\item $\langle \afa\bb\afa\bb\afa\BB\afa\BB \rangle$.
Apply \propc{cross corner}:
$$\langle \afa\bb\afa\bb\afa|\BB\afa\BB| \rangle \rightarrow
[\langle \afa\bb\afa\bb\afa| \rangle, \langle \BB\afa\BB| \rangle]
=[\langle \afa\bb\afa\bb \rangle, \langle \afa\BB \rangle]
= [\langle a^ib^ja^kb^l \rangle,\langle a^mB^n\rangle],$$
say, (so $\alpha(w) = i+k+m,~ \beta(w) = j+l+n$), and
$$\si(w) \le \si([\langle a^ib^ja^kb^l \rangle,\langle a^mB^n\rangle]) -1.$$
Now
$$\si([\langle a^ib^ja^kb^l \rangle,\langle a^mB^n\rangle])=
\si(\langle a^ib^ja^kb^l \rangle) + \si(\langle a^mB^n\rangle) +
\inn(\langle a^ib^ja^kb^l \rangle,\langle a^mB^n\rangle).$$
Because of the format of \coryc{lpcory}  we need to consider
two cases:
\begin{romlist}
\item  $i=k$ and $j=l=1$ (by the construction,
$i$ and $k$ cannot be 1). In that case $\si(\langle a^ib^ja^kb^l \rangle)
=1$. By  \propc{lp} (7), $\inn(\langle a^ib^ja^kb^l \rangle, \langle
a^mB^n\rangle) = (i+k)n + m(j+l)$, and by \propc{lp} (1)
$\si(\langle a^mB^n\rangle) = (m-1)(n-1)$. This gives
$$\si([\langle a^iba^ib\rangle,\langle a^mB^n\rangle]) =
1 + 2in + 2m + (m-1)(n-1)  = (2i-1)n +m (n+1) +2$$
and $\si(w) \leq (2i-1)n + m(n+1) +1$. On the other hand the word
$w' = \langle a^{2i}b^2a^mB^n\rangle$ has the same $\alpha$
and $\beta$-values as $w$
and (\propc{lp} (3))
$\si(w') = (2i+m-1)(n+1) = (2i-1)n +m(n+1) + (2i-1)$. Since as remarked 
above $i\geq 2$, it follows that $\si(w')\geq \si(w) + 2$.

\item For all other $\langle a^ib^ja^kb^l \rangle$, \coryc{lpcory} gives
$\si(\langle a^ib^ja^kb^l \rangle) \leq (i+k-1)(j+l-1)-2$, and
$$\si([\langle a^iba^ib\rangle,\langle a^mB^n\rangle]) 
\leq (i+k-1)(j+l-1)-2 + (i+k)n + m(j+l) + (m-1)(n-1) $$
$$= (i+k+m-1)(j+l+n-1)-1,$$
so $\si(w)\leq (i+k+m-1)(j+l+n-1)-2$. Comparing this estimate
with $\si(\langle a^{i+k}b^{j+l}a^mB^n\rangle) = (i+k+m-1)(j+l+n-1)$
(\propc{lp} (3) again) completes the proof.
\end{romlist}
\end{enumerate}
\end{proof}

\start{lem}{two} 
If $w$ uses exactly two letters and $w$
has three or more block-pairs, then there exists a word $w'$ with two 
fewer block-pairs, with the same $\alpha$ and
$\beta$ values, and such that $\si(w')\geq \si(w)+2$.
\end{lem}
\begin{proof} Suppose the two letters are $a$ and $b$, so
$w=\langle \afa\bb\afa\bb\afa\bb\dots\rangle$. Now proceed as in
the proof of \lemc{Three-letter 5b}.
\end{proof}

\subsection{End of the proof}
\start{prop}{upper bound h}
 Let $w$ be the reduced cyclic word corresponding to a
free homotopy class of curves on the punctured torus, 
with $h(w) = h > 0.$ 
Then there exists a word $w'$ such that 
$w'$ has one or two blocks,
$\alpha(w')=\alpha(w)$ and $\beta(w')=\beta(w)$, and  
$\si(w')~\ge~ \si(w)+h-2$.
\end{prop}

\begin{proof} 
If $h = 1 \mbox{~or~} 2$, then 
taking $w'=w$ satisfies the conclusions of the Proposition.

For $h>2$, we proceed by complete induction, and assume that the result 
holds for any word with a number of block-pairs smaller than $h$. 
Since $h$ is positive, $w$ contains 2, 3 or 4
distinct letters. We consider the cases separately.

\begin{numlist}
\item[(2 letters)] Suppose that $w$ contains exactly two distinct letters. 
If $w$ has more than two block-pairs, then by  
\lemc{two} there exists a word $v$ with 
$h-2$ block-pairs, with the same $\alpha$ and
$\beta$ values, and such that $\si(v)~\ge~ \si(w)+2$. 
By the  induction hypothesis, there exists $w'$ with one or two blocks, 
the same $\alpha$ and $\beta$ as $v$ and such that 
$\si(w') ~\ge~ \si(v)+(h-2)-2 ~\ge~ \si(w) +h-2,$ 
as desired. 

\item[(3 letters)] Suppose that $w$ contains exactly three distinct letters. 
If $h>4$,  the result follows from combining  \lemc{Three-letter 5b} 
and the induction hypothesis. If $h=4$ then
 $w$ must have one of the following forms:  
$\langle \afa\bb\afa\bb\afa\bb\afa\BB \rangle$, 
$\langle \afa\bb\afa\BB\afa\bb\afa\BB \rangle$ or  
$\langle \afa\bb\afa\bb\afa\BB\afa\BB \rangle$. 
\lemc{abababaB} 
covers these three cases. In the case $h=3$, the word can be supposed to be $w=\la a^{i}b^{j}a^{k}b^{l}a^{m}B^{n}\ra$. Taking $w'=\la a^{i+k}b^{j+l}a^{m}B^{n}\ra$, and applying \propc{lp}(3) and  (5), yields the desired result.

\item[(4 letters)] Suppose that $w$ contains all four letters  $a,b,A,B$. 
By \lemc{four}, there exists a word $v$ with $h-1$ block-pairs,  
with $\alpha(v)=\alpha(w),~\beta(v)=\beta(w)$ 
and such that $\si(v)~\ge~ \si(w)+1$. Now the result follows from our 
induction hypothesis. More explicitly, there exists a word $w'$, 
with same $\alpha$ and $\beta$ as $v$, with one or two block-pairs 
and such that $\si(w')~\ge~ \si(v)+(h-1)-2~\ge~ \si(w)+h-2$. 
\end{numlist}\end{proof}

\noindent
\begin{prooftext}{Proof of \theoc{main-th}.} 
If $h = 1 \mbox{~or~} 2$ the 
result follows from \propc{lp} (1-4). So suppose that $h > 2$,
By  \propc{upper bound h},  there exists a word $w'$ 
with $\alpha(w')=\alpha(w)$ and $\beta(w')=\beta(w)$,
$\si(w)\le \si(w')-h+2< \si(w')$, 
and such that $w'$ has 
one or two blocks. Referring to \propc{lp}(1-4), any such word
satisfies
$\si(w')\leq (\alpha(w')-1)(\beta(w')-1).$ This proves
part (1).

Part (2) of the theorem follows also, by inspection, from \propc{lp}(1-4).
\end{prooftext}

\subsection{Words with sub-maximal intersection number}

\start{lem}{ub4}  If $w$  is one of the following words: 
$\langle \mathsf{ababAB} \rangle,$ $\langle \mathsf{abAbaB} \rangle,$ $\langle \mathsf{abaBAB} \rangle,$ and $\langle \mathsf{abABaB} \rangle,$ then $$\si(w) \le  (\alpha(w)-1)(\beta(w)-1)-2.$$

\end{lem}

\begin{proof} 
\propc{lp} and \coryc{lpcory} can be applied 
after one or two cross-corner surgeries (\propc{cross corner}),
each of which increases the self-intersection number by at least one:

\begin{itemize}
\item  $\langle a^ib^ja^kb^l|A^mB^n| \rangle \rightarrow
\langle a^ib^ja^kb^l|b^na^m| \rangle = \langle a^{i+m}b^ja^kb^{l+n} 
\rangle$ 
\item
$\langle  a^ib^j|A^kb^la^mB^n| \rangle
\rightarrow
\langle  a^ib^j|b^nA^mB^la^k| \rangle 
=\langle  a^{i+k}b^{n+j}A^mB^l \rangle$

\item
$\langle a^ib^ja^k|B^lA^mB^n| \rangle \rightarrow
[\langle a^ib^ja^k|\rangle,\langle B^lA^mB^n|\rangle]
 =[\langle a^{i+k}b^j\rangle, \langle A^mB^{n+l} \rangle]$ 
\newline
$[\langle a^{i+k}|b^j\rangle, \langle A^m|B^{n+l} \rangle]
\rightarrow \langle a^{i+k}|a^mb^{n+l}|b^j\rangle =
\langle a^{i+k+m}b^{j+n+l}\rangle $. 
\item $ \langle a^i|b^jA^k|B^la^mB^n  \rangle \rightarrow
 \langle a^i|a^kB^j|B^la^mB^n  \rangle 
= \langle a^{i+k}B^{l+j}a^mB^n \rangle $ 
\end{itemize}
\end{proof}

The following lemma will be used in the proof of \theoc{upper bound 3}. 
Note that the special case it covers admits a 
bound for the self-intersection number sharper than that of 
\theoc{upper bound}.

\start{lem}{ub3} If $w$ is a word with three block-pairs, then $\si(w) \le (\alpha(w)-1)(\beta(w)-1)-2$. 
In particular,  if length $L =\alpha(w)+\beta(w)$ then
$$\si(w) \le \left \{ \begin{array}{ll} (L-2)^2/4-2 & \mbox{if $L$ is even,} \\
(L-1)(L-3)/4 -2&   \mbox{if $L$ is odd.}
\end{array} \right . $$ 
\end{lem}
\begin{proof} 
Without loss of generality,  we may suppose that the number $N$ of $\mathsf{A}$ and $\mathsf{B}$ blocks  in $w$ is at most three. 

If $N = 0$,  the result follows from 
\lemc{two} and \theoc{upper bound}.

If $N = 1$, we may suppose that $w= \langle  \mathsf{ababaB}  \rangle$ which  is covered by  \propc{lp}(5).

If $N = 2$,  we 
may
suppose that $w$ is one of $\langle \mathsf{abAbaB} \rangle$ or $\langle \mathsf{ababAB} \rangle$; 
if $N=3$, we 
may suppose  that $w$ is one of $\langle \mathsf{abaBAB} \rangle$ or
$\langle \mathsf{abABaB} \rangle$; for these cases, the result follows from \lemc{ub4}. 
\end{proof}

\start{lem}{abab bound}  Let $w$ be a word with two block-pairs
and two letters, say $a$ and $b$ with length $L \geq 4$.
Either $L=4$ and $w = abab$ with $\si(w)=1$, or
$$\si(w) \le \left \{ \begin{array}{ll} (L-2)^2/4-2 & \mbox{if $L$ is even,} \\
(L-1)(L-3)/4 -2&   \mbox{if $L$ is odd.}
\end{array} \right . $$ 
\end{lem}
\begin{proof} Refer to \coryc{lpcory}.

First notice that  $ab^jab^j$ has length $L=2j+2$, an even number,  and  if $j\ge 2$
$$(L-2)^2/4-2
= j^2-2 \geq 1 = \si(ab^jab^j).$$ So the lemma
holds for all words of the form $ab^jab^j$  and $a^{i}ba^{i}b$.

For the rest of the words in question,
$\si(w) \leq (i +k - 1)(j + l -1)-2$, so the result follows as in \remc{x1}.
\end{proof}

\start{theo}{upper bound 3} Let $w$ be a primitive reduced cyclic word 
of length $L>3$ and self-intersection number
$$\si(w) =\left \{ \begin{array}{ll} (L-2)^2/4-1 & \mbox{if $L$ is even,} \\
(L-1)(L-3)/4 -1&   \mbox{if $L$ is odd.}
\end{array} \right. $$ 
i.e. one less than the maximum possible for its length.
Then if $L$ is odd,   $w=r^{i}s^{j}R^{k}S^{l}$ with $i+k =\frac{L-1}{2}
\mbox{~or~} \frac{L+1}{2}$.

And if  $L$ is even,  $w$ has one of the following forms.
 \begin{numlist}
\item[{\rm (1)}] $\la r^{L/2-1}s^{L/2+1} \ra$; 
 \item[{\rm (2)}] $\la r^{i}s^{j}R^{k}S^{l} \ra$, $i+k= \frac{L}{2}$;
\item[{\rm (3)}] $\la r^{i}s^{j}r^{k}S^{l} \ra$, $i+k = \frac{L}{2}-1 \mbox{~or~}\frac{L}{2}+1.$
\end{numlist}
Here $r = a \mbox{~or~}A$ and $s =b \mbox{~or~}B$,
 or vice-versa.
\end{theo}

\start{rem}{a label} The primitive reduced cyclic words of length
$L\leq 3$, namely those of the form $a, ab, abb$ all have
self-intersection number zero, the maximum for those lengths
(cf. \tabc{counting table}). 
\end{rem}

\begin{prooftext}{Proof of \theoc{upper bound 3}} 
By \propc{upper bound h} and \lemc{ub3}, 
$h(w) = 1 \mbox{~or~}2$ (the only primitive words with zero block-pairs are
singletons, which do not satisfy the hypothesis.). 

We begin with the case $h(w)=2$. By \lemc{abab bound}, 
we can assume that $w = \langle \mathsf{abaB} \rangle \mbox{~or~} 
\langle \mathsf{abAB} \rangle$.

First suppose $w=\langle \mathsf{abaB} \rangle$. By \propc{lp} (3), 
$\si(w)=(\alpha(w)-1)(\beta(w)-1)$. 

If $L$ is odd, then  $\si(w)=(L-1)(L-3)/4-1$. 
Since $\beta(w)=L-\alpha(w)$, it follows that
$(L-1)(L-3)/4-1=(\alpha(w)-1)(L-\alpha(w)-1)$. 
This implies $\alpha(w)=(L\pm\sqrt{5})/2$, which is not an integer, 
a contradiction.

So $L$ is even, and    $(L/2-1)^{2}-1=(\alpha(w)-1)(L-\alpha(w)-1)$.  
This implies $\alpha(w) = \frac{n}{2}-1 \mbox{~or~} \frac{n}{2}+1$, as desired.

Now suppose $w=\langle \mathsf{abAB} \rangle$. 
The result follows from  \propc{lp} (4).
This settles the case $h(w)=2$.

If $h(w)=1$ then by \propc{lp}(1), 
$\si(w)=(\alpha(w)-1)(L-\alpha(w)-1)$.  The solutions of the equation  $$(\alpha(w)-1)(L-\alpha(w)-1)=(L-1)(L-3)/4-1$$ are 
$\alpha(w)=\frac{L-\sqrt{5}}{2}$ and  
$\alpha(w)=\frac{L+\sqrt{5}}{2}$. Hence there are no words of sub-maximal self-intersection with odd length $L$ and one block-pair.

On the other hand, the solutions of the equation 
$$(\alpha(w)-1)(L-\alpha(w)-1)=(L/2-1)^{2}-1$$ are $\alpha(w)=L/2-1$ and $L/2+1$; the result follows. 
\end{prooftext}

\start{theo}{counting 2} If $L$ is odd,  there are $(L-1)(L-3)$ 
distinct reduced cyclic words with self-intersection number
one less than the maximum for their length.

If $L$ is even, there are $5(L-2)^2/2$  distinct reduced cyclic words with self-intersection number
one less than the maximum for their length.
\end{theo}
\begin{proof} Refer to \theoc{upper bound 3}.
Suppose $L$ is odd. If $i+k=\frac{L-1}{2}$,  there are
$1, \dots, \frac{L-3}{2}$ possibilities for $i$, and  
$1 \dots, \frac{L-1}{2}$ possibilities for $j$. 
The total is $ \frac{L-3}{2} \frac{L-1}{2}$. Interchanging 
the roles of $i$ and $j$, and those of $a$ and $b$, 
we obtain $(L-1)(L-3)$.

Suppose $L$ is even:
there are $8$ words of the form $\la r^{L/2-1}s^{L/2+1} \ra$, together with
 $2(L/2-1)^{2}$ words of the form $\la r^{i}s^{j}R^{k}S^{l} \ra$ and 
$4L(L/2-2)$ words of the form $\la r^{i}s^{j}r^{k}S^{l} \ra$;
the total is  $5L^2/2-10L+10=5(L-2)^2/2$. 
\end{proof}

\start{rem}{xxxx} The leading coefficient of the polynomial expression for the number of maximal words of odd length is two times larger than that for even length, whereas for sub-maximal words the even leading coefficient is 2.5 times the odd leading coefficient. The discrepancies balance out to some extent, when one considers  maximal and sub-maximal words together.  For odd length $L$, this number is $3L^2-12L+17$, while for even length it is $7L^2/2-14L+18$. 
\end{rem}

\section{Computational results and conjectures for the pair of pants}

The ``pair of pants'' is the usual name for the surface with boundary,
homeomorphic to the thrice-punctured sphere. The same computational methods
that yielded \ctc{comp} suggest that the dependance of  
maximum self-intersection number on length for the pair of pants is
quadratic, just as it was for the punctured torus.

\start{ct}{pants1} For lengths $L\leq 15$, the 
maximal self-intersection
number of primitive reduced cyclic words of 
length $L$ on the pair of pants is:
$$\left \{ \begin{array}{ll} (L^{2}-1)/4 & \mbox{if $L$ is odd,} \\
L^{2}/4-1 &   \mbox{if $L \equiv 0 \pmod 4$,}\\
L^{2}/4-2 &   \mbox{if  $L>2$ and  $L \equiv 2 \pmod 4$,}\\
1 &   \mbox{if $L =2$.}\\
\end{array} \right . $$ 
Moreover, if $L$ is odd, the words realizing the maximal self-intersection number are 
 $r(rs)^\frac{L-1}{2}$, where $\{r,s\}=\{a, B\}$ or $\{r,s\}= \{A, b\}$.
(primitive words of even length follow a more complicated pattern, which cannot
be easily reduced to a formula).
\end{ct}

Removing the restriction ``primitive'' leads to:

\noindent
\start{ct}{pants2} For $L\leq 15$ 
the maximal self-intersection
number for a reduced cyclic word of 
length $L$ on the pair of pants is:
$$\left \{ \begin{array}{ll} (L^{2}-1)/4 & \mbox{if $L$ is odd,} \\
L^{2}/4+L/2-1 &   \mbox{if $L $ is even.}
\end{array} \right . $$ 
Moreover
\begin{numlist}
\item[{\rm (1)}] If $L$ is odd, the words realizing the maximal self-intersection number are 
 $r(rs)^\frac{L-1}{2}$, where $\{r,s\}=\{a, B\}$ or $\{r,s\}= \{A, b\}$
 \item[{\rm (2)}] If $L$ is even, the words realizing the maximal self-intersection number are 
 $(Ab)^\frac{L-1}{2}$ and  $(aB)^\frac{L-1}{2}$.
  \end{numlist}
\end{ct}

The next two computational theorems show radically different
behavior from what we know for the punctured torus. 

\noindent
\start{ct}{pants3} 
The number of distinct 
free homotopy classes of curves on 
the pair of pants of length $L$ realizing the
maximal self-intersection number is 
$$\left \{ \begin{array}{ll}2 & \mbox{if $L$ is even,} \\
4&   \mbox{if $L$ is odd.}
\end{array} \right . $$ 
\end{ct}

\start{ct}{pants4} For $L\leq 15$ the
{\em minimal} self-intersection number for the free homotopy
class on the pair of pants representing a primitive reduced
cyclic word of length $L$ is 0 for $L=1,2$ and $[L/2]$ (the
integer part of $L/2$) for $L\geq 3$. 
\end{ct}

It is reasonable to conjecture that all this behavior will persist for
higher values of $L$.

\start{rem}{other surfaces} Note that an analogue of \propc{cross corner} can be proved for any surface with boundary. So words with maximal self-interesection number cannot contain (the generalization of) diagonally opposed corners with reversed orientations. 
\end{rem}

\appendix
\section{Appendix: Proof of \propc{lp}}\label{appendix}

The seven parts of  \propc{lp} are proved separately as 
Propositions 
\ref{ab}, \ref{abab}, \ref{ab_aB}, \ref{abaB}, \ref{abAB}, 
\ref{ababaB}, \ref{abab_aB}.
The method of proof for each of these propositions is via
\theoc{clb si}: a counting of all occurrences
of each of the three types of linked pairs given in \defic{def}.

\begin{lplist}
\item These pairs are 
easy to count. They have the form $\{rr,ss\}$, where $r \in \{a,A\}$ and $s \in \{b,B\}$.
\item These have the form 
$\{{\it p}_1 \yy {\it p}_2, {\it q}_1  \yy {\it q}_2\}$, with $p_{1} \ne p_{2}$ and $q_{1}\ne q_{2}$. 
One locates all subwords $ \yy $ with 
two occurrences and ckecks for each pair if the
corresponding ${\it p}_1 \yy {\it p}_2$ and ${\it q}_1  \yy {\it q}_2$
are linked.
\item  
Analogously, these pairs are found by locating 
subwords  $ \yy $ which occur in our word or multiword along with their
inverse ${\mathsf Y}$. Such a pair will contribute to the count if the
corresponding ${\it p}_1\yy{\it p}_2$ and ${\it \overline{q}}_2\yy {\it \overline{q}}_1$ 
are linked, see \remc{xxx}.
\end{lplist}

\start{prop}{ab} $\si( \langle a^ib^j \rangle)=(i-1)(j-1)$.
\end{prop}
\begin{proof}  There are $i-1$ occurrences of $aa$ and $j-1$ occurrences
of $bb$ in  $\langle a^ib^j \rangle.$ Thus there are $(i-1)(j-1)$ pairs 
of type I. There are no pairs of the 
other two types.
\end{proof}

\start{prop}{abab}
$\si( \langle a^ib^ja^kb^l \rangle)=(i+k-2)(j+l-2)+ 1$ if $k=i$ and $j=l$; and 
$(i+k-2)(j+l-2)+|i-k|+|j-l|-1 $, otherwise.
\end{prop}

\begin{proof}
\vs

\begin{lplist}
\item There are $(i+k-2)(j+l-2)$  pairs of this kind.

\item In this case,  $ \yy \in \{a^K, B^K, a^Kb^J, b^Ka^J, a^Kb^Ja^L, b^Ja^Kb^L\}$ for some positive integers $J$, $K$ and $L$.
\begin{romlist}
\item$ \yy =a^K$.   Analysis:   \tabc{ababak}, using \tabc{tab-a}.  
The total number is  $|i-k|-1$ if $i\ne k$ and zero otherwise.
\begin{table}[htbp]
\centering
\begin{tabular}{|c|c|c|c|} \hline
 configuration&with&if&add\\\hline
$\{a^{k+2},ba^{k}b\}$ &$a^{k+2}$ in $a^i$& $k+2\le i$ & $i-k-1$\\\hline
$\{a^{i+2},ba^{i}b\}$ & $a^{i+2}$ in $a^k$& $i+2\le k$ & $k-i-1$\\\hline
\end{tabular}
\caption{Linked pairs in $\langle a^ib^ja^kb^l \rangle$ of type II with 
$ \yy =a^K$.}\label{ababak}

\end{table}

\item $ \yy = b^K$. With similar arguments as in the case (i), it can be shown that the number of pairs here is $|j-l|-1$ if $j\ne l$ and zero otherwise.

\item $ \yy =a^K  b^J$.  By \tabc{tab-ab}, the linked words of pairs with this $ \yy $ have the form $ba^Kb^Jb$ and 
$aa^Kb^Ja$. Analysis: \tabc{akbj}.

\begin{table}[htbp]
\centering
\begin{tabular}{|c|c|c|c|} \hline
 configuration&with&if&add\\\hline
 $\{ba^{i}b^{l}b,aa^{i}b^{l}a\}$&$ba^{i}b^{l}b$ in $ba^{i}b^{j}$, $aa^{i}b^{l}a$ in $a^{k}b^{l}a$ &$i < k$ and $j>l$& 1\\\hline 
$\{ba^{k}b^{j}b,aa^{k}b^{j}a\}$& $ba^{k}b^{j}b$ in $ba^{k}b^{l}$, $aa^{k}b^{j}a$ in $a^{i}b^{j}a$ &$k<i$ and $j<l$&1 \\\hline
\end{tabular} 
\caption{Linked pairs in $\langle a^ib^ja^kb^l \rangle$ of type II with 
$\yy =a^K  b^J$.The three types of linked pairs can be added as follows:}\label{akbj}
\end{table}

\item $ \yy =$$b^Ka^J$.   By \tabc{tab-ab}, the linked words have the form $ab^Ka^Ja$  and 
$bb^Ka^Jb$. Analysis: \tabc{bkaj}.
\begin{table}[htbp]
\centering
\begin{tabular}{|c|c|c|c|} \hline
 configuration&with&if&add\\\hline
$\{bb^{j}a^{i}b,ab^{j}a^{i}a\}$& $bb^{j}a^{i}b$ in $b^{l}a^{i}b$,  $ab^{j}a^{i}a$ in  $ab^{j}a^{k}$ &$i < k$ and $j<l$& 1\\\hline
 $\{bb^{l}a^{k}b,ab^{l}a^{k}a\}$& $bb^{l}a^{k}b$ in $b^{j}a^{k}b$, $ab^{l}a^{k}a$ in $ab^{l}a^{i}$ & $k<i$ and $j>l$& 1\\\hline
\end{tabular}
\caption{Linked pairs in $\langle a^ib^ja^kb^l \rangle$ of type II  with 
$\yy =b^Ka^J$.}\label{bkaj}
\end{table}

By (iii) and (iv) we add $1$ if $k\ne i$ and $j\ne l$.

\item $ \yy =$$a^Kb^Ja^L$. Since $\yy$ has two occurrences, $j=l$. 
By  \tabc{tab-aba} the linked pairs have the form 
$\{aa^{K}b^{j}a^{L}a,ba^{K}b^{j}a^{L}b\}$. 
There are two pairs of this kind, namely  $\{aa^{i}b^{j}a^{k}a, 
ba^{i}b^{j}a^{k}b\}$ and $\{aa^{k}b^{l}a^{i}a, ba^{k}b^{l}a^{i}b\}$. 
Each of the possibilities  implies that $i<k$ and $k<i$. 
Hence, such pairs are not possible.

\item $ \yy  = {\it b^Ka^Jb^L}$. 
As in case (v), there are no linked pairs of this form.
\end{romlist}

\item  There are no pairs of type III  because the
word contains no occurrence of a letter and its inverse.

\end{lplist}

If $i=k$ and $j=l$ add $1$,  because the word has the form $w^{2}$, where $w$ is a primitive word. 
Adding up all the  contributions completes the proof. 
\end{proof}

\start{prop}{ab_aB}
$\inn( \langle a^ib^j \rangle,  \langle a^m\B^n \rangle)=in+mj$.
\end{prop}
\begin{proof}
\vs

 \begin{lplist}
\item  There are $(i-1)(n-1)+(j-1)(m-1)$ linked pairs of this type.
\item   $\yy = a^K$ for some positive integer $K$. Analysis:
\tabc{aak}, using \tabc{tab-a}.
The contributions of the different rows may be grouped in the following way: 
(a + c + d)$=i-1$, (b + f + h)$=m-1$ and 
(e + g + i)$ =1$.
\begin{table}[htp]\centering
\begin{tabular}{|c|c|c|c|c|} \hline
 &configuration&with&if&add\\\hline
a& $\{a^{K}b, Ba^{K}\}$& $a^{K}b$ in $a^{i}b$, $Ba^{K}$ in $Ba^{m}$ &$K\in\{2\dots \min(m,i)\} $&
$\min(m,i)-1$  \\\hline
b& $\{a^{i}B,  ba^{i}\}$&  $ba^{K}$ in $ba^{i}$, $a^{K}B$ in $a^{m}B$&$K\in\{2\dots \min(m,i)\} $ &
 $\min(m,i)-1$  \\\hline
c& $\{a^{m+2}, Ba^{m}B\}$& $a^{m+2}$ in $a^i$&  $m+2 \le i$& $i-m-1$\\\hline
d& $\{a^{m+1}b,Ba^{m}B\}$& $a^{m+1}$ in $a^i$ & $m+1 \le i$ & $1$\\\hline
e& $\{ba^{m+1},Ba^{m}B\}$& $a^{m+1}$ in $a^i$&  $m<i$ & $1$\\\hline
f& $\{a^{i+1}B,ba^{i}b\}$& $a^{i+1}$ in $a^m$&  $i<m$ & $1$\\\hline
g& $\{Ba^{i+1},ba^{i}b\}$& $a^{i+1}$ in $a^m$&  $i<m$ & $1$\\\hline
h& $\{a^{i+2},ba^{i}b\}$& $a^{i+2}$ in $a^m$&   $i+2\le m$& $m-i-1$\\\hline
i& $\{ba^{i}b,Ba^{m}B\}$&   & $i=m$ & $1$\\\hline
\end{tabular}
\caption{Linked pairs of  $\langle a^ib^j \rangle$ and   $\langle a^m\B^n \rangle$ of type II with $\yy=a^{K}$.}\label{aak}
\end{table}

\item   $ \yy = b^K$. Combining \remc{xxx} with \tabc{tab-a} we  
analyze these pairs in \tabc{aA}. Here $(\mbox{a + c + f})=n-1, 
(\mbox{b + g + i})=j-1$ and $(\mbox{d + e + h})= 1$.
\begin{table}[htp]\centering
\begin{tabular}{|c|c|c|c|c|} \hline
 &configuration&with&if&add\\\hline
a& $\{ab^K, aB^K\}$     & $ab^K$ in $ab^j$, $aB^K$ in $aB^n$&$K\in\{2\dots \min(j,n)\} $ &   $\min(j,n)-1$\\\hline
b& $\{b^Ka, B^Ka\}$     & $ab^K$ in $ab^j$, $aB^K$ in $aB^n$&$K\in\{2\dots \min(j,n)\} $  &   $\min(j,n)-1$\\\hline
c&  $\{ab^ja, B^{j+1}a\}$ &$B^{j+1}a$ in $B^na$           & $j<n$ & $1$\\\hline 
d& $\{ab^ja, aB^{j+1}\}$&   $aB^{j+1}$ in $aB^n$          &$j<n$ &  $1$\\\hline 
e& $\{ab^ja, aB^na\}$   &                                 &$j=n$ & $1$\\\hline
f& $\{ab^ja, B^{j+2}\}$ &  $B^{j+2}$ in $B^{n}$                           & $j+2 \le n$ & $n-j-1$\\\hline
g& $\{b^{n+1}a, aB^na\}$&          $b^{n+1}a$ in $b^{j}a$                   &$n<j $  &  $1$\\\hline
h& $\{ab^{n+1}, aB^na\}$&          $ab^{n+1}$ in $ab^{j}$             &$n<j $&   $1$\\\hline
i& $\{b^{n+2}, aB^na\}$ & $b^{n+2}$ in $b^{j} $                      &$n+2 \le j $ & $j-n-1$ \\\hline
\end{tabular} 
\caption{Linked pairs of  $\langle a^ib^j \rangle$ and   $\langle a^m\B^n \rangle$ of type III with $\yy=b^{K}$.}\label{aA}
\end{table}
 \end{lplist}
 Adding the contributions from each of the three types yields the result.
\end{proof}

\start{prop}{abaB} $\si( \langle a^ib^ja^k\B^l \rangle)=(i+k-1)(j+l-1)$.
\end{prop}
\begin{proof}
\vs

\begin{lplist}
\item  There are 
$(i+k-2)(j-1)$
ocurrences of pairs $\{aa,bb\}$ and $(i+k-2)(l-1)$
ocurrences of pairs $\{aa,BB\}$. The total is
$(i+k-2)(j+l-1)$.

\item  $ \yy  \in \{a^I, b^J,B^K\}$ for some positive 
integers $I,J$ and $K$. There are no linked pairs with $ \yy  \in \{ a^I, b^J,B^K\}$. By \tabc{tab-a}, (interchanging roles of $a$'s and $b$'s) there are no pairs such that $\yy \in \{b^{J},B^{K}\}$.  
We analyze each of the possible pairs with $\yy=a^{I}$ in Table \ref{abaB1},
using \tabc{tab-a}. We show that there are $i+k-2$ of this type.
\begin{table}[htp]
\centering
\begin{tabular}{|c|c|c|c|c|} \hline
 &configuration&with&if&add\\\hline
& $\{a^{I}b,Ba^{I}\}$& $a^{I}b$ in $a^{i}b$, $Ba^{I}$ in $Ba^{i}$ &$I \in \{2,3,\dots, i\}$ & $i-1$\\\hline
& $\{ba^{I},a^{I}B\}$& $ba^{I}$ in $ba^{i}$, $a^{I}B$ in $a^{i}B$ & $I\in  \{2,3,\dots, k\}$ & $k-1$\\\hline
\end{tabular}
\caption{Linked pairs of   $\langle a^ib^ja^k\B^l \rangle$ of  type II 
with $ \yy =a^{I}$.}\label{abaB1}
\end{table}

\item   $ \yy =b^J$. Analysis: \tabc{b^i}. 
The total contribution from the Type III linked pairs comes to $j+l-1$. 
\begin{table}[htp] \centering
\begin{tabular}{|c|c|c|c|c|}
 \hline
 &configuration&with&if&add\\\hline
 & $\{b^{J}a,B^{J}a\}$ & $b^{J}a$ in $b^{j}a$, $B^{J}a$ in $B^{l}a$ &  $J \in \{2,3,\dots\min(j,l)\}$ &    $\min(j,l)-1$\\\hline 
 &  $\{ab^{J},aB^{J}\} $&$ab^{J}$ in $ab^{j}$, $aB^{J}$ in $aB^{l}$ & $J\in \{2,3,\dots \min(j,l)\}$&  $\min(j,l)-1$\\\hline
 & $\{a b^{j}a,B^{j+1}a\} $& $B^{j+1}a$ in $B^{l}a$ & $j<l$&   1\\\hline
 & $\{a b^{j}a,aB^{j+1}\} $& $aB^{j+1}$ in $aB^{l}$& $j<l$ & 1 \\\hline
 & $\{a b^{j}a,B^{j+2}\} $&$B^{j+2}$ in $B^{l}$ & $j+2 \le l$ &  $l-j-1$ \\\hline
 & $\{a B^{l}a,b^{l+1}a\} $& $b^{l+1}a$ in $b^{j}$ &$j>l$ &1\\\hline 
 & $\{a B^{l}a,ab^{l+1}\} $&$ab^{l+1}$ in $b^{j}$ &$j>l$ &1\\\hline 
 & $\{a B^{l}a,b^{l+2}\} $& $b^{l+2}$ in $b^{j}$ &$j>l+1$& $j-l-1$\\\hline
 & $\{a B^{l}a,ab^{j}a\} $& &$j=l$ &1\\\hline
\end{tabular}
\caption{Linked pairs of  $\langle a^ib^ja^k\B^l \rangle$ of  type III 
with $\yy=b^{J}$.} 
\label{b^i}
 \end{table}
 
\end{lplist}
It follows that the
total self-intersection number of $ \langle a^ib^ja^kB^l  \rangle$ is
$(i+k-2)(j+l-2)+i+k-2+j+l-1 = (i+k-1)(j+l-1)$. 
\end{proof}

\start{prop}{abAB} $\si(\langle a^{i}b^{j}A^{k}B^{l} \rangle=(i+k-1)(j+l-1)-1$.
\end{prop}
\begin{proof}
\vs

 \begin{lplist}
\item   There are in total
$(i+k-2)(j+l-2)$ pairs of this type.
\item In this case $ \yy  \in \{a^{I}, b^{J},B^{K},A^L\}$ for some integers 
$I, J, K$ and $L$. Analysis: \tabc{abABx}, using \tabc{tab-a}. 

\begin{table}[htp] \centering
\begin{tabular}{|c|c|c|c|c|}
 \hline
 &configuration&with&if&add\\\hline
 &$\{a^{I}b,Ba^{I}\}$& $a^{I}$ in $a^{i}$ & $I\in \{2,3,\dots, i\}$&$i-1$\\\hline 
 &$\{bA^{I},A^{I}B\}$&$A^{I}$ in $A^{k}$ &$I\in \{2,3,\dots, k\} $&$k-1$\\\hline
 &$\{b^{I}A,ab^{I}\}$& $b^{I}$ in $b^{j}$& $I\in \{2,3,\dots, j\}$&$j-1$\\\hline
 &$\{AB^{I},B^{I}a\}$&$B^{I}$ in $B^{l}$ &$I\in \{2,3,\dots, l\}$ &$l-1$\\\hline
\end{tabular}
\caption {Linked pairs of $\langle a^{i}b^{j}A^{k}B^{l} \rangle$ of type II.}
\label{abABx}\end{table}

\item $\yy \in\{a^i, b^i, A^i,B^i\}$ since the
inverses  of the subwords  $\{ab, bA, AB, Ba\}$ are $\{BA, aB, ba, Ab\}$,
which do not occur in $\langle a^{p}b^{q}A^{r}B^{s} \rangle$. 
It follows from \remc{xxx} and \tabc{tab-a} that there are no linked pairs
of this type.
\end{lplist}
Total of I, II, III:
$$
(i+k-2)(j+l-2)-i+j+k+l-4=(i+k-1)(j+l-1)-1
$$
\end{proof}

\start{prop}{ababaB}
$$\si( \langle a^ib^ja^kb^la^mB^n \rangle)=(i+k+m-1)(j+l+n-1)
-2    (k+\min(j,l)-1).$$
\end{prop}
\begin{proof}

\begin{lplist}
\item    There are  $(i+k+m-3)$ ocurrences of  $aa$, $(j+l-2)$ ocurrences of $bb$  and  $(n-1)$
 occurrences of $BB$. This gives 
$(i+k+m-3)(j+l+n-3)$ linked pairs of  type I.
\item  $\yy \in \{a^Ib^J, b^Ja^I,a^Ib^Ja^K, b^{J},B^{I},a^{I}\}$ for some integers $I,J$ and $K$. By \remc{corners} we do not need to consider words of the form 
$rs {\mathsf u}sR$. 

 \begin{romlist}
\item $\yy=a^{I}b^J$. We analyze these linked pairs in \tabc{ab7}. The total contribution of these pairs is $1$ if $l>j$.
\begin{table}[htp] \centering
\begin{tabular}{|c|c|c|c|c|}
 \hline
 &configuration&with&if&add\\\hline
 & $\{aa^{I}b^{J}a,Ba^{I}b^{J}b\}$  & & &not linked\\\hline
 & $\{aa^{I}b^{j}a,ba^{k}b^{J}b\}$ & $aa^{I}b^{j}a$ in $a^{i}b^{j}a$, $ba^{k}b^{J}b$ in $ba^{k}b^{l}$ & $l>j$  and  $i>k$  & 1\\\hline
 & $\{aa^{I}b^{J}b,Ba^{i}b^{j}a\}$ & $aa^{I}b^{J}b$ in $a^{k}b^{l}$ & $l>j$ and $i<k$  &1\\\hline
 & $\{Ba^{i}b^{j}a,ba^{k}b^{J}b\}$ & $ba^{k}b^{J}b$ in $ba^{k}b^{l}$ & $l>j$ and $i=k$ & 1\\\hline
 & $\{ba^{I}b^{J}a,aa^{I}b^{J}b\}$ & & &not linked\\\hline
 & $\{ba^{I}b^{J}a,Ba^{i}b^{J}b\}$ & & &not linked\\\hline
\end{tabular}
\caption{Linked  pairs of $\langle a^ib^ja^kb^la^mB^n \rangle$ of type II with $\yy=a^{I}b^J$.}\label{ab7}
\end{table}

\item $\yy= b^Ja^I$. We analyze these pairs in \tabc{ab8}. The total here is $1$ if $l<j$.
 \begin{table}[htp] \centering
\begin{tabular}{|c|c|c|c|c|}
 \hline &configuration&with&if&add\\\hline
 & $\{ab^{l}a^Ia,bb^{J}a^{k}b\}$ & $ab^{l}a^{I}a$ in $ab^{l}a^{m}$,   $bb^{J}a^{k}b$ in $b^{j}a^{k}b$ & $l<j$  and $k<m$ &  1\\\hline
 & $\{ab^{J}a^Ia,bb^{J}a^{m}B\}$ &  & &   not linked\\\hline
 & $\{ab^{J}a^Ib,bb^{J}a^{I}a\}$ & & &   not linked\\\hline
 & $\{ab^{J}a^Ib,bb^{J}a^{m}B\}$ & & &   not linked\\\hline
 & $\{ab^{l}a^mB,bb^{J}a^{k}b\}$ &$bb^{J}a^{k}b$ in $b^{j}a^{k}b$ & $l<j$ and $m=k$& 1\\\hline
 & $\{ab^{l}a^mB,bb^{J}a^{I}a\}$ & $bb^{J}a^{I}a$ in $b^{j}a^{k}$ & $l<j$ and $m<k$& 1\\\hline
\end{tabular}
\caption{Linked  pairs of $\langle a^ib^ja^kb^la^mB^n \rangle$ of type II with $\yy= b^Ja^I$.}\label{ab8}
\end{table}

\item $\yy= a^{I}b^{J}a^{K}$. These pairs are analyzed in \tabc{ab9} where it is shown that there are no pairs of this type.
\begin{table}[htp] \centering
\begin{tabular}{|c|c|c|c|c|}
 \hline
 &configuration&with&if&add\\\hline
& $\{Ba^{i}b^{j}a^{I}, ba^{k}b^{l}a^{m}B\}$ & & & 
 not linked (\remc{corners})\\\hline
& $\{Ba^{i}b^{j}a^{I}, a^{I}b^{l}a^{m}B\}$ & & &   not linked\\\hline
& $\{Ba^{i}b^{j}a^{k}b, a^{I}b^{l}a^{m}B\}$ & & & 
not linked (\remc{corners})\\\hline
& $\{a^{I}b^{j}a^{I}, ba^{k}b^{l}a^{m}B\}$ & & &
not linked (\remc{corners})\\\hline
& $\{a^{I}b^{j}a^{k}b, ba^{k}b^{l}a^{m}\}$ & & &   not linked   \\\hline
\end{tabular}
\caption{Linked  pairs of $\langle a^ib^ja^kb^la^mB^n \rangle$ of type II with $\yy= a^{I}b^{J}a^{K}$.}\label{ab9}
\end{table}
 
\item $\yy=b^{I}$. These pairs are analyzed in \tabc{ab10}. They 
contribute $|l-j|-1$ if $|l-j|>1$.
\begin{table}[htp] \centering
\begin{tabular}{|c|c|c|c|c|}
 \hline
 &configuration&with&if&add\\\hline
 & $\{b^{I},ab^{j}a\}$ & $b^{I}$ in $b^{l}$ & $j+1<l$ & $l-j-1$\\\hline
 & $\{b^{I},ab^{l}a\}$ & $b^{I}$ in $b^{j}$  & $l+1<j$ & $j-l-1$\\\hline
 & $\{ab^{I},b^{I}a\}$ & & & not linked\\\hline
\end{tabular}
\caption{Linked  pairs of $\langle a^ib^ja^kb^la^mB^n \rangle$ of type II with $\yy=b^{I}$,}\label{ab10}
\end{table}

\item $\yy=B^{I}$. The only possible pair, $\{aB^{I},B^{I}a\}$   is not linked.
\item  $\yy=a^{I}$. Here (see \tabc{tab-a}) the candidates for linked 
subwords are $ \{a^{I}b,ba^{I},Ba^{I},a^{I}B,a^{I},ba^{I}b\}.$
The number of occurrences of each 
possible pair is tabulated in \tabc{ab11}. The total number is $2m+2i-4$.
\begin{table}[htp] \centering
\begin{tabular}{|c|c|c|c|c|}
 \hline
 &configuration&with&if&add\\\hline
& $\{a^{I}b,Ba^{I}\}$& $Ba^I$ in $Ba^i$, $a^Ib$ in $a^ib$&$I \in \{2,3,\dots,i\}$ &$i-1$\\\hline
&$\{a^{I}b,Ba^{I}\}$ &$Ba^I$ in $Ba^i$, $a^Ib$ in $a^kb$&$I \in \{2,3,\dots,\min(i,k)\}$ &$\min(k,i)-1$\\\hline
&$\{ba^{I},a^{I}B\}$  &$a^IB$ in $a^mB$, $ba^I$ in $ba^m$& $I \in \{2,3,\dots,m\}$& $m-1$ \\\hline 
& $\{ba^{I},a^{I}B\}$ &$a^IB$ in $a^mB$, $ba^I$ in $ba^k$&$I \in \{2,3,\dots,\min(k,m)\}$ & $\min(k,m)-1$\\\hline 
& $\{aa^{k}B,ba^{k}b\}$& $aa^{k}B$ in $a^{m}B$ &$m>k$& 1\\\hline
& $\{Ba^{k+1},ba^{k}b\}$& &$i>k$& 1\\\hline
&$\{ba^{k}b,a^{k+2}\}$ & $a^{k+2}$ in $a^i$&$i >k+1$&$i-k-1$\\\hline
&$\{ba^{k}b,a^{k+2}\}$ & $a^{k+2}$ in $a^m$&$m>k+1$&$m-k-1$\\\hline
\end{tabular}
\caption{Linked  pairs of $\langle a^ib^ja^kb^la^mB^n \rangle$ of type II 
with $\yy=a^{I}$.}\label{ab11}
\end{table}
\end{romlist}

The total number of linked pairs of type II is 
$|j-l|+2m+2i-4$.

\item $\yy=b^{I}$. Analysis: \tabc{ab14}, using \remc{xxx} and \tabc{tab-a}.
 The values in \tabc{ab14} may be grouped as follows: (a+b+e+f+i+j)=$2n-2$, 
(d+k+q)=$j-1$, (c+l+r)=$l-1$, (g+h+m+n+o+p)=$2$.

\begin{table}[htp] \centering
\begin{tabular}{|c|c|c|c|c|}
 \hline
 &configuration&with&if&add\\\hline
a & $\{B^{K}a, b^Ka\}$ &$B^{K}$ in $B^{n}$, $b^{K}$ in $b^{l}$ & $K  \in \{2,3,\dots, \min(n,l)\}$& $\min(n,l)-1$\\\hline
b & $\{B^{K}a, b^Ka\}$ & $B^{K}$ in $B^{n}$, $b^{K}$ in $b^{j}$&  $K \in \{2,3,\dots, \min(n,j)\}$ &$\min(n,j)-1$\\\hline
c & $\{aB^{K}, ab^K\}$ & $B^{K}$ in $B^{n}$, $b^{K}$ in $b^{l}$&$K\in \{2,3,\dots, \min(n,l)\}$ & $\min(n,l)-1$\\\hline
d & $\{aB^{K}, ab^K\}$ & $B^{K}$ in $B^{n}$, $b^{K}$ in $b^{j}$& $K \in \{2,3,\dots, \min(n,j)\}$  & $\min(n,j)-1$\\\hline
e & $\{B^{K}a, ab^ja\}$&$B^{K}$ in $B^{n}$ &$n>j$&1\\\hline
f & $\{B^{K}a, ab^la\}$&$B^{K}$ in $B^{n}$ &$n>l$&1\\\hline
g & $\{aB^{K}, ab^ja\}$& $B^{K}$ in $B^{n}$ &$n>j$&1\\\hline
h & $\{aB^{K}, ab^la\}$& $B^{K}$ in $B^{n}$ &$n>l$&1\\\hline
i & $\{B^{K}, ab^ja\}$&$B^{K}$ in $B^{n}$ & $n>j+1$ & $n-j-1$\\\hline
j & $\{B^{K}, ab^la\}$&$B^{K}$ in $B^{n}$ & $n>l+1$ & $n-l-1$\\\hline
k&$\{aB^{n}a, bb^Ka\}$&  $bb^{K}$ in $b^{j}$&$j>n$&1\\\hline
l &$\{aB^{n}a, bb^Ka\}$ & $bb^{K}$ in $b^{l}$&$l>n$&1\\\hline
m&$\{aB^{n}a, ab^Kb\}$&$b^{K}b$ in $b^{j}$ &$j>n$&1\\\hline
n &$\{aB^{n}a, ab^Kb\}$ &$b^{K}b$ in $b^{l}$ &$l>n$&1\\\hline
o&$\{aB^{n}a, ab^ja\}$& &$n=j$&1\\\hline
p &$\{aB^{n}a, ab^la\}$ & & $n=l$&1\\\hline
q&$\{aB^{n}a, b^K\}$& $b^{K}$ in $b^{j}$ &$j>n+1$&$j-n-1$\\\hline
r &$\{aB^{n}a, b^K\}$ & $b^{K}$ in $b^{l}$&$l>n+1$&$l-n-1$\\\hline
\end{tabular}
\caption{Linked  pairs of $\langle a^ib^ja^kb^la^mB^n \rangle$ of type III
 with $\yy=b^{K}$.}\label{ab14}
\end{table}
\end{lplist}

Total for types I, II, III:
$$(i+k+m-3)(j+l+n-3)+|j-l|+2m+2i-4+2n+j+l-2. $$
Applying $|j-l|=\max(j,l)-\min(j,l)$ and $j+l = \max(j,l)+\min(j,l)$ 
yields the desired result.

\end{proof}

\start{prop}{abab_aB}
$\inn( \langle a^ib^ja^kb^l\rangle, \langle a^m\B^n \rangle)=(i+k)n+m(j+l)$
\end{prop}
\begin{proof}
\vs

 \begin{lplist}
\item  There are $(i+k-2)(n-1)+(m-1)(j+l-2)$ of these.

\item  $ \yy = a^K$ for some positive integer $K$. 
Analysis: \tabc{ab16}, using \tabc{tab-a}.
The values in \tabc{ab16} may be grouped as follows:
(a+e+g)=$i-1$,
 (b+f+h)=$k-1$,
 (c+d+k+l+o+p)= $2m-2$, and 
 (i+j+m+n+p+q)= $2$
 
\begin{table}[htp]
\centering
\begin{tabular}{|c|c|c|c|c|} \hline
 &configuration&with&if&add\\\hline

a &$\{a^{K}b,Ba^{K}\}$  &$a^Kb$ in $a^ib$&$K \in \{2,3,\dots, \min(m,i)\}$ & $\min(m,i)-1$\\\hline
b& $\{a^{K}b,Ba^{K}\}$& $a^Kb$ in $a^kb$&$K \in \{2,3,\dots, \min(m,k)\}$  & $\min(m,k)-1$\\\hline
 c&$\{a^{i}B,ba^{i}\}$ & $ba^K$ in $ba^i$&$K \in \{2,3,\dots, \min(m,i)\}$  & $\min(m,i)-1$\\\hline
d&$\{a^{i}B,ba^{i}\}$ & $ba^K$ in $ba^k$& $K \in \{2,3,\dots, \min(m,k)\}$ & $\min(m,k)-1$\\\hline
e&  $\{a^{m+2},Ba^{m}B\}$& $a^{m+2}$ in $a^i$&  $m+2 \le i$ & $i-m-1$\\\hline
f & $\{a^{m+2},Ba^{m}B\}$ & $a^{m+2}$ in $a^k$& $m+2 \le k$ & $k-m-1$\\\hline
g&$\{a^{m+1}b,Ba^{m}B\}$ & $a^{m+1}b$ in $a^ib$&   $m<i$ & $1$\\\hline
h&$\{a^{m+1}b,Ba^{m}B\}$ & $a^{m+1}b$ in $a^jb$&   $m<k$ & $1$\\\hline

i& $\{ba^{m+1},Ba^{m}B\}$& $ba^{m+1}$ in $ba^i$& $m<i$& $1$\\\hline
j &$\{ba^{m+1},Ba^{m}B\}$ & $ba^{m+1}$ in $ba^j$& $m<k$& $1$\\\hline

k &$\{a^{k+1}B, ba^{k}b\}$ & $a^{k+1}B$ in $a^m$ & $k<m$ & 1 \\\hline 

l & $\{a^{i+1}B,ba^{i}b\}$ & $a^{i+1}B$ in $a^mB$ &$i<m$ & $1$\\\hline
m&$\{Ba^{k+1},ba^{k}b\}$& $B^{k+1}a$ in $Ba^m$ &$k<m$ & $1$ \\\hline
n&$\{Ba^{i+1},ba^{i}b\}$& $B^{i+1}a$ in $Ba^m$ &$k<m$ & $1$ \\\hline

o& $\{a^{i+2},ba^{i}b\}$ & $a^{i+2}$ in $a^m$&  $i+2 \le m$ & $m-i-1$\\\hline
 p& $\{a^{i+2},ba^{k}b\}$ & $a^{i+2}$ in $a^m$ &   $k+2 \le m$ & $m-k-1$\\\hline

q& $\{ba^{i}b,Ba^{m}B\}$& & $i=m$ & $1$\\\hline
r& $\{ba^{k}b,Ba^{m}B\}$ & & $k=m$ & $1$\\\hline
\end{tabular}
\caption{Linked  pairs of  $\langle a^ib^ja^kb^l\rangle$ and  
$\langle a^m\B^n \rangle$ of type II  with $ \yy = a^K$.}\label{ab16}
\end{table}

\item  $\yy=b^K$. Analysis:  \tabc{ab15}, using  \tabc{tab-a} (interchanging the roles of $a$'s and $b$'s) and \remc{xxx}. The values group in the following way: 
 (a+b+e+f+k+l)= $2n-2$, 
(c+m+q)= $j-1$,
(d+p+r)= $l-1$,
(g+h+i+j+n+o) = $2$.
\begin{table}[htp] \centering
\begin{tabular}{|c|c|c|c|c|} \hline
 &configuration&with&if&add\\\hline
a& $\{ab^K, aB^K\}$& $b^{K}$ in $b^j$& $K\in \{2,3,\dots,\min(j,n)\}$ & $\min(j,n)-1$\\\hline
b& $\{ab^K, aB^K\}$& $b^{K}$ in $b^l$&  $K\in \{2,3,\dots,\min(l,n)\}$  &  $\min(l,n)-1$\\\hline
c&$\{b^Ka, B^Ka\}$ &  $b^{K}a$ in $b^ja$&   $K\in \{2,3,\dots,\min(j,n)\}$ &   $\min(j,n)-1$\\\hline
d& $\{b^Ka, B^Ka\}$& $b^{K}a$ in $b^la$ &  $K\in \{2,3,\dots,\min(l,n)\}$  &   $\min(l,n)-1$\\\hline
e&$\{ab^{l}a,B^{l+1}a\}$& &$l< n$ & $1$\\\hline 
f & $\{ab^{j}a,B^{j+1}a\}$& &$j<  n$&  $1$\\\hline 
g & $\{ab^{j}a,aB^{j+1}\}$& &$j<  n$&  $1$\\\hline 
h &$\{ab^{l}a,aB^{l+1}\}$& &$l<  n$ & $1$ \\\hline 
i & $\{ab^ja, aB^na\}$& &$j=n$ &  $1$\\\hline 
j &$\{ab^{l}a,aB^{n}a\}$& &$l=n$ & $1$\\\hline 
k &  $\{ab^ja, B^{j+2}\}$& &$j+2 \le n$ & $n-j-1$ \\\hline 
l &$\{ab^{l}a,B^{l+2}\}$& &$l+2 \le n$ & $n-l-1$ \\\hline 
m & $\{b^{n}a,aB^{n}a\}$& $b^na$ in $b^ja$&$n<j$& $1$\\\hline
n & $\{b^{n}a, aB^{n}a\}$& $b^na$ in $b^la$& $n<l$ & $1$ \\\hline
o & $\{ab^K, aB^Ka\}$& &$n<j$& 1\\\hline 
p &    $\{ab^K, aB^Ka\}$               & &$n<l$ &1\\\hline 
q& $\{b^{n+2}, aB^na\}$& $b^{n+2}a$ in $b^{j}$&$n+2 \le j$ &$j-n-1$ \\\hline 
r &$ \{b^{n+2}, aB^na\}$   &$b^{n+2}a$ in $b^{l}$ &$n+2 \le l$ &  $l-n-1$ \\\hline
\end{tabular} 
\caption{Linked  pairs of   $\langle a^ib^ja^kb^l\rangle$ and  
$\langle a^m\B^n \rangle$  of type III with $\yy=b^K$.}\label{ab15}
\end{table}

\end{lplist}

Grand Total: $(i+k-2)(n-1)+(m-1)(j+l-2)+i+j+k+l+2m+2n-4 = (i+k)n +m(j+l) $.  
\end{proof}

%


\bibliographystyle{amsalpha}

\textsc{
Department of Mathematics,
Stony Brook University,
Stony Brook, NY, 11794.
}

\emph{E-mail address}{:\;\;}\texttt{moira@math.sunysb.edu, tony@math.sunysb.edu}
\end{document}